\documentclass[11pt]{article}
\usepackage{amsfonts,amsmath,subfigure}
\usepackage{amssymb,mathrsfs,latexsym,bm}
\usepackage{arydshln}

\begin{document}
\newcommand{\qed}{\hphantom{.}\hfill $\Box$\medbreak}
\newcommand{\proof}{\noindent{\bf Proof \ }}
\newtheorem{Theorem}{Theorem}[section]
\newtheorem{Lemma}[Theorem]{Lemma}
\newtheorem{Corollary}[Theorem]{Corollary}
\newtheorem{Remark}[Theorem]{Remark}
\newtheorem{Example}[Theorem]{Example}
\newtheorem{Definition}[Theorem]{Definition}
\newtheorem{Construction}[Theorem]{Construction}

\thispagestyle{empty}
 \renewcommand{\thefootnote}{\fnsymbol{footnote}}

\begin{center}
{\Large\bf Two-dimensional balanced sampling plans avoiding
adjacent units \footnote{Supported by the Fundamental Research Funds for the Central Universities grant $2013$JB$Z005$, Shandong Excellent Young Scientist Research
Award Fund Project grant BS2012SF005, the Natural Science Foundation of Ningbo grant $2013$A$610102$, the K. C. Wong Magna Fund at Ningbo University, the NSFC  grants $11271042$, $11201252$ and $61373007$}}

\vskip12pt
Xiaomiao Wang$^1$, Tao Feng$^2$, Jing Zhang$^3$ and  Yan Xu$^1$ \\[2ex] {\footnotesize
$^1$Department of Mathematics, Ningbo University, Ningbo 315211, P. R. China
\\$^2$Department of Mathematics, Beijing Jiaotong University, Beijing 100044, P. R. China
\\$^3$College of Mathematics and Physics, Qingdao University of Science and Technology, Qingdao 266042, P. R. China}\\
{\footnotesize
 wangxiaomiao@nbu.edu.cn, tfeng@bjtu.edu.cn, zhangjing@qust.edu.cn}
\vskip12pt

\end{center}

\vskip12pt

\noindent {\bf Abstract:} Hedayat et al. first introduced balanced sampling plans for the exclusion of contiguous units. Wright detailed the results of a preliminary investigation of two-dimensional balanced sampling plans avoiding
adjacent units (2-BSAs), and pointed out explicitly three types of 2-BSAs, which have different adjacency scheme, namely ``Row and Column'', ``Sharing a Border'' and ``Island''. This paper will provide more details for the three types of 2-BSAs from the point of view of design theory.



\section{Introduction}

In environmental and ecological populations, neighboring units within a finite population, spatially or sequentially ordered, may provide similar information. It is intuitively appealing to select a sample that avoids the selection of adjacent units. Balanced sampling plans excluding adjacent units have been proposed as a means of achieving such a goal.

When the units are arranged in a one-dimensional ordering, the population may follow a circular ordering, in which the first unit of the population is contiguous
with the last unit, or a linear ordering, in which the first
unit is not contiguous with the last unit. Hedayat et al. (1988a,b) first proposed a sampling plan for a given
circular population of size $N$, for which a sample size $k$ is obtained without replacement such that the second-order inclusion probabilities are $0$ for contiguous units and some positive constant for non-contiguous units. Stufken (1993) extended this to
balanced sampling avoiding adjacent units.

Suppose the population set $X$ is identified with $Z_N:= \{0,1,\ldots,N-1\}$. Let $m$ be
a positive integer. Two units $x$ and $y$ are said to be {\em adjacent} if
$x-y\in\{-m,-m+1,\ldots,-1, 0, 1,\ldots,m\}$, where the arithmetic is performed modulo $N$. A {\em circular one-dimensional balanced sampling plan avoiding adjacent units}, or simply a 1-BSA$(N,k,\lambda;m)$, is a sampling plan of $k$-subsets (called {\em blocks}) from the population of size $N$ such that two units that are adjacent do not appear together in any block while any two non-adjacent units appear together in exactly $\lambda$ blocks. When $m=1$, it is often referred to as a {\em one-dimensional balanced sampling plan excluding contiguous units}, denoted by 1-BSEC$(N,k,\lambda)$.

Hedayat et al. (1988a) observed that if a 1-BSEC$(N,k,\lambda)$ exists, then so does a 1-BSEC$(N+3,k,\lambda')$. Using this observation with several small values of $N$, they showed that for $k=3$ or $4$ and $N\geq 3k$, there exists a 1-BSEC$(N,k,\lambda)$ for some $\lambda$. Colbourn and Ling (1998, 1999) gave a complete existence theorem for $k=3$ and $4$. We only quote the following result for the later use.

\begin{Theorem}{\rm (Colbourn and Ling 1998)}\label{one-dimensional BSEC}
A $1$-BSEC$(N,3,\lambda)$ exists if and only if $N\geq 9$ and $\lambda (N-3)\equiv 0\ ({\rm mod}\ 6)$.
\end{Theorem}

A 1-BSA$(N,k,\lambda;m)$ is equivalent to a special case of a partial balanced incomplete block design, called polygonal design. In terms of polygonal designs, Stufken and Wright (2001) proved that for $k\in\{5,6,7\}$ and $N\geq 3k+1$, there exists a 1-BSEC$(N,k,\lambda)$ for some $\lambda$, with the possible exception of $N=22$ and $k=7$. For more information on balanced sampling plans with small $m$ and $k$, the reader may refer to Wright and Stufken (2008), Iqbal et al. (2009). We remark that Wright and Stufken (2008) also gave a discussion systematically on linear one-dimensional balanced sampling plan avoiding adjacent units; we will not provide any detail here.

In this paper, we focus on two-dimensional balanced sampling plans. Two-dimensional populations will be restricted to those consisting of $r$ rows and $c$ columns, and units within such populations will be identified as ordered pairs $\{(i,j): 0\leq i \leq r-1; 0\leq j\leq c-1\}$. While the general concept of adjacency under one-dimensional populations is easily extended to two-dimensional populations, enormous flexibility is gained in the application of the concept. Wright (2008) proposed three possible adjacency schemes that can be considered under two-dimensional populations in Fig. $1$, where the units labeled $\heartsuit$ for a given adjacency scheme are considered to be adjacent to the $\spadesuit$ unit. A two-dimensional balanced sampling plan will be denoted by 2-BSA$(r,c,k,\lambda;{\rm adjacency\ scheme})$, where $k$ is the sample size, and $\lambda$ is the number of samples containing two given non-adjacent units.

\begin{figure}[t]
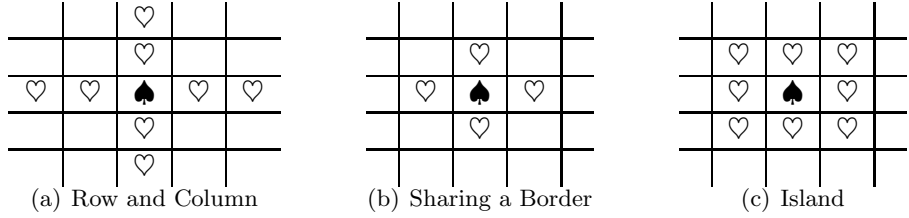
\renewcommand\figurename{Fig.}
\centering
\subtable[Row and Column]{
\begin{tabular}{l|l|l|l|l}
&&$\heartsuit$&&\\\hline
&&$\heartsuit$&&\\\hline
$\heartsuit$&$\heartsuit$&$\spadesuit$&$\heartsuit$&$\heartsuit$\\\hline
&&$\heartsuit$&&\\\hline
&&$\heartsuit$&&\\
\end{tabular}
}
\qquad
\subtable[Sharing a Border]{
\begin{tabular}{l|l|l|l|l}
&&&&\\\hline
&&$\heartsuit$&&\\\hline
&$\heartsuit$&$\spadesuit$&$\heartsuit$&\\\hline
&&$\heartsuit$&&\\\hline
&&&&\\
\end{tabular}
}
\qquad
\subtable[Island]{
\begin{tabular}{l|l|l|l|l}
&&&&\\\hline
&$\heartsuit$&$\heartsuit$&$\heartsuit$&\\\hline
&$\heartsuit$&$\spadesuit$&$\heartsuit$&\\\hline
&$\heartsuit$&$\heartsuit$&$\heartsuit$&\\\hline
&&&&\\
\end{tabular}
}
\caption{Graphical representations of adjacency schemes for two-dimensional populations}
\end{figure}

We shall details the three kinds of two-dimensional balanced sampling plans in Sections 2, 3 and 4, respectively.

\section{2-BSAs with adjacency scheme ``Row and Column''}

First we give a clear mathematical definition for 2-BSAs with adjacency scheme ``Row and Column''. Given $(x,y)\in Z_r\times Z_c$, the points $(i,y)$ and $(x,j)$ for any $i\in Z_r$ and $j\in Z_c$ are said to be {\em row-column-mates} of the point $(x, y)$.

A {\em two-dimensional balanced sampling plan avoiding row-column-mates} is a pair $(X,{\cal B})$, where $X=Z_r\times Z_c$ and $\cal B$ is a collection of $k$-subsets of $X$ (called {\em blocks}) such that any two points that are row-column-mates do not appear in any block while any two points that are not row-column-mates appear in exactly $\lambda$ blocks. It is simply denoted by a 2-BSA$(r,c,k,\lambda;{\rm RC})$.

From the point of view of design theory, it is closely related to a special kind of design, called modified group divisible design. Let $r,c,k$ and $\lambda$ be positive integers. A {\em modified group divisible design} (MGDD) is a quadruple $(X,{\cal G},{\cal H},{\cal B})$ which satisfies the following properties:
\begin{enumerate}
    \item[(1)] $X$ is a finite set of $rc$ {\em points};
    \item[(2)] ${\cal G}$ is a partition of $X$ into $r$ subsets (called {\em groups}), each of size $c$;
    \item[(3)] ${\cal H}$ is another partition of $X$ into $c$ subsets (called {\em holes}), each of size $r$, such that $|H\cap G|=1$ for each $H\in {\cal H}$ and $G\in {\cal G}$;
    \item[(4)] $\cal B$ is a set of subsets (called {\em blocks}) of $X$, each of
    size $k$, such that no block contains two distinct points of any group or any hole, but any other pair of distinct points of $X$ occurs in exactly $\lambda$ block of $\cal B$.
\end{enumerate}
Such a design is denoted by a $(k,\lambda)$-MGDD of type $c^r$. Obviously it is nothing but a 2-BSA$(r,c,k,\lambda;{\rm RC})$.

\begin{Example}\label{3-MGDD(3^3)}
Here we give an example of a $(3,1)$-MGDD of type $(3^4)$. We label the $12$ points as follows:
\begin{center}
\begin{tabular}{lll}
$(0,0)$&$(0,1)$&$(0,2)$\\
$(1,0)$&$(1,1)$&$(1,2)$\\
$(2,0)$&$(2,1)$&$(2,2)$\\
$(3,0)$&$(3,1)$&$(3,2)$\\
\end{tabular}
\end{center}
Each row is a group and each column is a hole. The $12$ blocks are
\begin{center}
\begin{tabular}{ll}
$\{(0,i),(1,i+1),(2,i+2)\}$,&
$\{(0,i),(1,i+2),(3,i+1)\}$,\\
$\{(0,i),(2,i+1),(3,i+2)\}$,&
$\{(1,i),(2,i+2),(3,i+1)\}$.\\
\end{tabular}
\end{center}
where $i=0,1,2$, and the arithmetic is reduced modulo $3$. It is also a $2$-BSA$(4,3,3,1;{\rm RC})$.
\end{Example}

Assaf (1990) first introduced the notion of MGDDs and settled the existence when $k=3$. The existence of $(4,\lambda)$-MGDDs was investigated by Assaf (1997), Assaf and Wei (1999), Ling and Colbourn (2000), Ge et al. (2003). We summarize their results in the language of two-dimensional balanced sampling plans as follows.

\begin{Theorem} \label{MGDD}
\begin{enumerate}
\item[$(1)$]
There exists a $2$-BSA$(r,c,3,\lambda;{\rm RC})$ if and only if $r,c\geq 3$, $\lambda(r-1)(c-1)\equiv 0\ ({\rm mod}\ 2)$, and $\lambda r(r-1)c(c-1)\equiv 0\ ({\rm mod}\ 3)$.
\item[$(2)$] There exists a $2$-BSA$(r,c,4,\lambda;{\rm RC})$ if and only if $r,c\geq 4$, $\lambda(r-1)(c-1)\equiv 0\ ({\rm mod}\ 3)$, except when $\lambda=1$ and $\{r,c\}=\{6,4\}$.
\end{enumerate}
\end{Theorem}

We remark that there are also some results on $(5,\lambda)$-MGDDs; see Abel and Assaf (2002, 2008).

\section{2-BSAs with adjacency scheme ``Sharing a Border''}

Bryant et al. (2002) first detailed a study on 2-BSAs with adjacency scheme ``Sharing a Border''. For $(x,y)\in Z_r\times Z_c$, the points
$(x-1,y)$, $(x+1,y)$, $(x,y-1)$, and $(x,y+1)$ (reducing the arithmetics modulo $r$ and $c$ in the first and second coordinates, respectively) are said to be {\em $2$-contiguous} to the point $(x, y)$.

A {\em two-dimensional balanced sampling plan avoiding $2$-contiguous units} is a pair $(X,{\cal B})$, where $X=Z_r\times Z_c$ and $\cal B$ is a collection of $k$-subsets of $X$ (called {\em blocks}) such that any two $2$-contiguous points do not appear in any block while any two points that are not $2$-contiguous appear in exactly $\lambda$ blocks. It is denoted by a 2-BSA$(r,c,k,\lambda;{\rm SB})$, or simply a $2$-BSEC$(r,c,k,\lambda)$ as used by Bryant et al. (2002).

When $r=1$ or $c=1$, a 2-BSEC can be seen as a balanced sampling plan excluding contiguous units with only one dimension. In this case each point has only two $2$-contiguous points. So a 2-BSEC$(1,c,k,\lambda)$ is just a 1-BSEC$(c,k,\lambda)$. If we allow $r$ or $c$ to be 2, then no point would have four 2-contiguous points, and such design has little significance in applications, so we always assume that $r,c\geq 3$ in this section unless otherwise specified.

Bryant et al. (2002) established the necessary and sufficient conditions for the existence of $2$-BSEC$(r,c,3,1)$s. By using modified group divisible designs, they gave an infinite family for the existence of $k=3$ and general $\lambda$. Much less is known about the existence of $2$-BSEC$(r,c,4,\lambda)$s; see Ge et al. (2003), Kong et al. (2008).

\begin{Theorem}\label{3-BSEC}{\rm (Bryant et al. 2002)}
\begin{enumerate}
\item[$(1)$] Let $r,c\geq 3$. There exists a $2$-BSEC$(r,c,3,1)$ if and only if $r$ and $c$ are odd and either $r\equiv c\equiv 3\ ({\rm mod}\ 6)$ or $r\not\equiv c\ ({\rm mod}\ 6)$.
\item[$(2)$] Let $r,c\in\{3\}\cup\{i:i\geq 9\}$. If $\lambda(r-3)\equiv 0\ ({\rm mod}\ 6)$ and $\lambda(c-3)\equiv 0\ ({\rm mod}\ 6)$, then there exists a $2$-BSEC$(r,c,3,\lambda)$.
\end{enumerate}
\end{Theorem}


In this section, we present the necessary and sufficient conditions for the existence of $2$-BSEC$(r,c,3,\lambda)$s. In the sampling context, the particular value of $\lambda$ may be less important than the ease with which blocks can be selected. However, designs with smaller $\lambda$ require less storage to represent explicitly for this purpose. By counting the number of blocks and the number of blocks containing a given point in a $2$-BSEC$(r,c,3,\lambda)$, the following necessary condition for its existence can be obtained.

\begin{Lemma}{\rm (Bryant et al. 2002)} \label{necessity for any index}
Let $r,c\geq3$. If a $2$-BSEC$(r,c,3,\lambda)$ exists, then $\lambda rc(rc-5)\equiv 0\ ({\rm mod}\ 6)$, and $\lambda (rc-5)\equiv 0\ ({\rm mod}\ 2)$.
\end{Lemma}

\begin{Lemma}\label{nonexistence}
No $2$-BSEC$(4,3,3,\lambda)$ exists for any $\lambda$.
\end{Lemma}

\proof Assume that there were a $2$-BSEC$(4,3,3,\lambda)$ on $Z_4\times Z_3$. Then the blocks should be the four types via the first coordinate: $\{(0,a_1),(1,a_2),(2,a_3)\}$, $\{(0,b_1),(1,b_2),(3,b_3)\}$, $\{(0,c_1),(2,c_2),(3,c_3)\}$, $\{(1,d_1),(2,d_2),(3,d_3)\}$. Let $x_1$, $x_2$, $x_3$ and $x_4$ be the number of blocks of the four types, respectively. So $x_1+x_2+x_3+x_4=14\lambda$. On the other hand, observing that the total number of the pairs corresponding to the first coordinate $(0,2)$ is $9\lambda$, so $x_1+x_3=9\lambda$. Similarly, by calculating the number of pairs corresponding to the first coordinate $(1,3)$, we have $x_2+x_4=9\lambda$. Thus  $x_1+x_2+x_3+x_4=18\lambda$, a contradiction. \qed

We shall show that the necessary conditions for the existence of a 2-BSEC$(r,c,3,\lambda)$ are also sufficient except when $(r,c)=(4,3)$. The union of a 2-BSEC$(r,c,3,\lambda_1)$ and a 2-BSEC$(r,c,3,\lambda_2)$ is a 2-BSEC$(r,c,3,\lambda_1+\lambda_2)$. Hence it suffices to establish the existence of 2-BSECs for the minimum value of $\lambda$. Note that by the symmetry of $r$ and $c$, there is no difference between a 2-BSEC$(r,c,3,\lambda)$ and a 2-BSEC$(c,r,3,\lambda)$ essentially.

\subsection{Combinatorial tools}

\subsubsection{Difference method}

Difference method plays an important role in the direct construction for designs. The distinguishing feature of this method is that the properties of a design can easily be obtained from the sets of shifts instead of listing all blocks of the design.

Suppose $B\subset Z_r\times Z_c$. Let $\Delta B=\{(a_1-a'_1,a_2-a'_2):(a_1,a_2),(a'_1,a'_2)\in B,(a_1,a_2)\neq(a'_1,a'_2)\}$, where the arithmetic is reduced modulo $r$ and $c$ in the first and second coordinates, respectively. Let $\lambda$ be a positive integer and $S$ a set. We denote by $\lambda\cdot S$ a multiset containing each element of $S$ exactly $\lambda$ times. The following lemma is simple but very useful; its proof is straightforward and thus omitted here.

\begin{Lemma}\label{difference method}
Suppose there exist $k$-subsets $B_1,B_2,\ldots,B_b$ of $Z_r\times Z_c$ such that
 $$\bigcup_{i=1}^b \Delta B_i=\lambda\cdot(Z_r\times Z_c\setminus\{(0,0),(0,\pm1),(\pm1,0)\}). $$
 Then there exists a $2$-BSEC$(r,c,k,\lambda)$.
\end{Lemma}

The subsets $B_1,B_2,\ldots,B_b$ of $Z_r\times Z_c$ in Lemma \ref{difference method} are called {\em base blocks} of the 2-BSEC. In this paper, to save space, for each element $(x,y)$ of $Z_r\times Z_c$, we sometimes simply write $x_y$ instead of $(x,y)$.

\begin{Example}\label{11,4,3,2}
There exists a $2$-BSEC$(11,4,3,2)$. Only base blocks are listed below:
\vspace{0.1cm}\\{\tabcolsep 0.05in
\begin{tabular}{lllllll}
$\{0_0,0_2,5_3\}$, & $\{0_0,1_1,4_2\}$, & $\{0_0,1_2,4_0\}$,&
$\{0_0,1_3,4_2\}$, & $\{0_0,1_1,4_0\}$, & $\{0_0,1_2,5_1\}$,&
$\{0_0,1_3,5_2\}$, \\ $\{0_0,2_0,4_1\}$, & $\{0_0,2_2,4_1\}$,&
$\{0_0,2_0,5_0\}$, & $\{0_0,2_1,5_3\}$, & $\{0_0,2_2,5_2\}$,&
$\{0_0,2_3,5_0\}$.
\end{tabular}}
\end{Example}

\subsubsection{Holey group divisible design}

Holey group divisible designs can be thought of as a natural generalization of modified group divisible designs. Let $g_1,g_2,\ldots,g_t$ and $n$ be positive integers. A {\em holey group divisible design} (HGDD) is a quadruple $(X,{\cal G},{\cal H},{\cal B})$ which satisfies the following properties:
\begin{enumerate}
    \item[(1)] $X$ is a finite set of $n\sum_{i=1}^t g_i$ {\em points};
    \item[(2)] ${\cal G}=\{G_1,G_2,\ldots,G_n\}$ is a partition of $X$ into $n$ subsets, (called {\em groups}), each of size $\sum_{i=1}^t g_i$;
    \item[(3)] ${\cal H}=\{H_1,H_2,\ldots,H_t\}$ is another partition of $X$ into $t$ subsets, (called {\em holes}), such that $|H_i|=ng_i$ and $|H_i\cap G_j|=g_i$ for any $1\leq j\leq n$;
    \item[(4)] $\cal B$ is a set of $k$-subsets (called {\em blocks}) of $X$, such that no block contains two distinct points of any group or any hole, but any other pair of distinct points of $X$ occurs in exactly $\lambda$ blocks of $\cal B$.
\end{enumerate}
If $\cal H$ contains $u_l$ holes of size $nh_l$, $1\leq l\leq s$,
then we call $(n,h_1^{u_1} h_2^{u_2}\cdots h_s^{u_s})$ the {\em type} of the HGDD. Such a design is denoted by a $(k,\lambda)$-HGDD of type $(n,h_1^{u_1} h_2^{u_2}\cdots h_s^{u_s})$. A $(k,\lambda)$-HGDD of type $(n,1^m)$ is just a $(k,\lambda)$-MGDD of type $m^n$.

\begin{Theorem}\label{3-HGDD-wei}{\rm (Wei 1993)}
There exists a $(3,\lambda)$-HGDD of type $(n,h^u)$ if and only if $n,u\geq 3$, $\lambda(u-1)(n-1)h\equiv 0\ ({\rm mod}\ 2)$ and $\lambda u(u-1)n(n-1)h^2\equiv 0\ ({\rm mod}\ 3)$.
\end{Theorem}

\begin{Theorem}\label{3-HGDD}{\rm (Wang and Yin 1999)}
There exists a $(3,1)$-HGDD of type $(n,h^u w^1)$ if and only if $(1)$ $n\geq 3$, $u=2$ and $h=w$; or $(2)$ $n\geq 3$, $u\geq 3$, $0\leq w\leq h(u-1)$, $hu(n-1)\equiv 0\ ({\rm mod}\ 2)$, $(n-1)(w-h)\equiv 0\ ({\rm mod}\ 2)$ and $hun(n-1)(h(u-1)-w)\equiv 0\ ({\rm mod}\ 3)$.
\end{Theorem}

To present our construction for 2-BSECs via HGDDs, we introduce a new configuration called quasi-modified group divisible designs. Let $r$ and $c$ be positive integers. Given $(x,y)\in Z_r\times Z_c$, the points $(x+1,y)$, $(x-1,y)$ (reducing the arithmetics modulo $r$) and $(x,j)$ for any $j\in Z_c$ are said to be {\em related} to the point $(x, y)$. We define a new configuration $(X,{\cal B})$, where $X=Z_r\times Z_c$ and $\cal B$ is a collection of $k$-subsets of $X$ (called {\em blocks}) such that any two related points do not appear in any block while any two points that are not related appear in exactly $\lambda$ blocks. Such a configuration is called a {\em quasi-modified group divisible design}, and denoted by a $(k,\lambda)$-QMGDD of type $c^r$. The sets $\{i\}\times Z_c$, $i\in Z_r$, are called the {\em groups} of the QMGDD. When $c=2,3$, a $(k,\lambda)$-QMGDD of type $c^r$ is just a 2-BSEC$(r,c,k,\lambda)$.

Similar to Lemma \ref{difference method}, the following result is straightforward.

\begin{Lemma}\label{difference method-QMGDD}
Let $E=\{(1,0),(-1,0)\}\cup\{(0,j):j\in Z_c\}$. Suppose there exist $k$-subsets $B_1,B_2,\ldots,B_b$ of $Z_r\times Z_c$ such that
 $$\bigcup_{i=1}^b \Delta B_i=\lambda\cdot(Z_r\times Z_c\setminus E).$$
Then there exists a $(k,\lambda)$-QMGDD of type $c^r$.
\end{Lemma}

The subsets $B_1,B_2,\ldots,B_b$ of $Z_r\times Z_c$ in Lemma \ref{difference method-QMGDD} are called {\em base blocks} of the QMGDD.  In this paper, the notation $\alpha\{a_1,a_2,a_3\}$ always means $\alpha$ copies of the block $\{a_1,a_2,a_3\}$.

\begin{Example}\label{example QMGDD-(2^5)}
There exists a $(3,\lambda)$-QMGDD of type $c^r$ for $(r,c,\lambda)\in\{(5,2,1),(8,5,2),(7,7$, $3),(7,2,6)\}$. Only base blocks are listed below.
\vspace{0.1cm}\\
\begin{tabular}{llllll}
$(r,c,\lambda)=(5,2,1):$&
$\{0_0,1_1,3_0\}$.\vspace{0.1cm}\\
$(r,c,\lambda)=(8,5,2):$
\end{tabular}
\\
\begin{tabular}{llllll}
$\{0_0,4_0,7_1\}$, & $\{0_0,2_0,6_2\}$, &
$\{0_0,2_0,3_2\}$, & $\{0_0,3_0,5_2\}$,&
$\{0_0,3_0,4_2\}$, & $\{0_0,3_1,7_2\}$,\\
$\{0_0,1_1,6_2\}$, & $\{0_0,2_1,7_2\}$,&
$\{0_0,1_1,3_2\}$, & $\{0_0,4_1,2_2\}$,&
$\{0_0,6_1,5_2\}$.\vspace{0.1cm}\\
\end{tabular}
\\\hspace*{0.2cm}$(r,c,\lambda)=(7,7,3):$
\\
\begin{tabular}{llllll}
$\{0_0,0_2,0_4\}$, & $\{0_0,0_3,1_4\}$, &
$\{0_0,0_3,2_4\}$, & $\{0_0,0_2,3_6\}$,&
$\{0_0,1_1,3_3\}$, & $\{0_0,1_1,3_2\}$,\\
$\{0_0,1_2,3_5\}$, & $\{0_0,1_2,3_4\}$,&
$\{0_0,1_3,3_4\}$, & $2\{0_0,1_3,3_6\}$,&
$\{0_0,1_5,3_2\}$, & $2\{0_0,1_5,3_3\}$, \\
$\{0_0,1_2,3_1\}$, & $\{0_0,1_4,3_2\}$,&
$\{0_0,1_4,3_1\}$, & $\{0_0,1_6,3_1\}$,&
$2\{0_0,1_6,3_5\}$.\vspace{0.1cm}\\
\end{tabular}
\\\hspace*{0.2cm}$(r,c,\lambda)=(7,2,6):$
\\
\begin{tabular}{llllll}
$2\{0_0,1_1,3_0\}$,& $2\{0_0,1_1,3_1\}$,& $2\{0_0,1_1,4_0\}$, &$2\{0_0,2_0,4_1\}$, $\{0_0,2_0,4_0\}$,& $\{0_0,2_1,4_0\}$.\\
\end{tabular}
\end{Example}

\begin{Lemma}\label{small QMGDD}
There exists a $(3,1)$-QMGDD of type $2^8$.
\end{Lemma}

\proof All $32$ blocks can be obtained by developing the following $16$ blocks by $(-,+1\ {\rm mod}\ 2)$:
\begin{center}
\begin{tabular}{llllll}
$\{0_0,2_0,4_0\}$,&$\{1_0,3_0,5_0\}$,&$\{0_0,3_0,4_1\}$,&
$\{0_0,5_0,7_1\}$,&$\{0_0,6_0,1_1\}$,&$\{1_0,4_0,5_1\}$,\\
$\{1_0,6_0,3_1\}$,&$\{1_0,7_0,4_1\}$,&$\{2_0,5_0,0_1\}$,&
$\{2_0,6_0,5_1\}$,&$\{2_0,7_0,1_1\}$,&$\{3_0,6_0,0_1\}$,\\
$\{3_0,7_0,2_1\}$,&$\{4_0,6_0,2_1\}$,&$\{4_0,7_0,6_1\}$,&
$\{5_0,7_0,3_1\}$.
\end{tabular}
\end{center}

\begin{Construction}\label{construction from hgdd}
Suppose that there exists a $(k,\lambda)$-HGDD of type $(n,h_1^{u_1} h_2^{u_2}\cdots h_s^{u_s})$. If there exist a $(k,\lambda)$-QMGDD of type $h_i^n$ for each $1\leq i\leq s$, and a $1$-BSEC$(\sum_{i=1}^s h_iu_i,k,\lambda)$, then there exists a $2$-BSEC$(n,\sum_{i=1}^s h_iu_i,k,\lambda)$.
\end{Construction}

\proof Let $(X,\cal{G},\cal{H},\cal{B})$ be the given $(k,\lambda)$-HGDD of type $(n,h_1^{u_1} h_2^{u_2}\cdots h_s^{u_s})$. For each $G\in{\cal G}$, $|G|=\sum_{i=1}^s h_iu_i$, and we construct a $1$-BSEC$(|G|,k,\lambda)$ on the set $G$. Denote the set of its blocks by ${\cal D}_G$.

For $H\in\cal{H}$, according to the definition  of HGDD, $|H\cap G|=h_i$ for any $G\in{\cal G}$ and some $1\leq i\leq s$. Now for each $H\in\cal{H}$, we construct a $(k,\lambda)$-QMGDD of type $h_i^n$ on the set $H$ with groups $H\cap G$, $G\in{\cal G}$. Denote the set of its blocks by ${\cal A}_H$.

Let ${\cal C}=(\bigcup_{H\in{\cal H}}{\cal A}_H)\cup(\bigcup_{G\in{\cal G}}{\cal D}_G)\cup{\cal B}$.
It is readily checked that $(X,{\cal C})$ is the required $2$-BSEC$(\sum_{i=1}^s h_iu_i,n,k,\lambda)$.  \qed

\subsubsection{Incomplete group divisible design}

An {\em incomplete group divisible design} (IGDD) is a quadruple $(X,Y,{\cal G},{\cal B})$, where $X$ is a set of points, $Y$ is a subset of $X$ (called the {\em hole}), $\cal G$ is a partition of $X$ into {\em groups}, and $\cal B$ is a collection of subsets of $X$ (called {\em blocks}) such that
\begin{enumerate}
\item[$(1)$] for each block $B\in{\cal B}$, $|B\cap Y|\leq 1$;
\item[$(2)$] no pair of points of $Y$ occurs in any block;
\item[$(3)$] any pair of points from $X$ which are not both in $Y$ occurs either in same group or in exactly $\lambda$ blocks, but not both.
\end{enumerate}
A $(k,\lambda)$-IGDD of type $(v_1,h_1)^{u_1}(v_2,h_2)^{u_2}\cdots(v_s,h_s)^{u_s}$ is an IGDD in which every block has size of $k$ and there are $u_i$ groups of size $v_i$, each of which intersects the hole in $h_i$ points, $i=1,2,\ldots,s$. When $Y=\emptyset$, an incomplete group divisible design is often called a {\em group divisible design}, and we use the notation $(k,\lambda)$-GDD of type $v_1^{u_1}v_2^{u_2}\cdots v_s^{u_s}$ instead of $(k,\lambda)$-IGDD of type $(v_1,0)^{u_1}(v_2,0)^{u_2}\cdots (v_s,0)^{u_s}$.

\begin{Theorem}{\rm (Zhu 1993)}\label{3-GDD}
There exists a $(3,\lambda)$-IGDD of type $(m,h)^u$ if and only if $m\geq 2h$, $\lambda m(u-1)\equiv 0\ ({\rm mod}\ 2)$, $\lambda (m-h)(u-1)\equiv 0\ ({\rm mod}\ 2)$ and $\lambda u(u-1)(m^2-h^2)\equiv 0\ ({\rm mod}\ 6)$.
\end{Theorem}


\begin{Theorem}{\rm (Heinrich and Zhu 1986)}\label{4-IGDD}
For $v\geq 3h$ and $h\geq 1$, there is a $(4,1)$-IGDD of type $(v,h)^4$ except when $v=6$ and $h=1$.
\end{Theorem}

As a straightforward corollary of Theorem 4.2 in Kong et al. (2008), we have the following construction.

\begin{Construction}\label{Construction from IGDD}
Suppose there exists a $(k,\lambda)$-IGDD of type $(v_1,2)^{u_1}(v_2,2)^{u_2}\cdots(v_s,2)^{u_s}$. Let $u=\sum_{i=1}^s u_i$. If there exist a $(k,\lambda)$-GDD of type $2^u$, a $(k,1)$-MGDD of type $k^r$, and a $2$-BSEC$(r,v_i,k,\lambda)$ for each $1\leq i\leq s$, then there exists a $2$-BSEC$(r,\sum_{i=1}^s v_iu_i,k,\lambda)$.
\end{Construction}

\subsection{$\lambda=2$}

\begin{Lemma}\label{3c+x recur}
Let $r\geq 3$, $c\geq 6$ and $3\leq x\leq c$. If there exist a $2$-BSEC$(r,c,3,2)$ and a $2$-BSEC$(r,x,3,2)$, then there exists a $2$-BSEC$(r,3c+x,3,2)$.
\end{Lemma}

\proof For $c\geq 6$, take a $(4,1)$-IGDD of type $(c,2)^4$ from Theorem \ref{4-IGDD}. Truncate one group to $x$ points and it is required that the two points in the hole are not removed. Now replace each block not containing truncated points by the blocks of a $(3,2)$-GDD of type $1^4$. Replace each block containing truncated points by two copies of the block. This yields a $(3,2)$-IGDD of type $(c,2)^3(x,2)^1$. Start from this IGDD and apply Construction \ref{Construction from IGDD} to obtain a 2-BSEC$(r,3c+x,3,2)$, where the needed $(3,2)$-GDD of type $2^4$ is from Theorem \ref{3-GDD} and the needed $(3,1)$-MGDD of type $3^r$ is from Theorem \ref{MGDD}(1). \qed

\begin{Lemma}\label{m=3 index 2}
There exists a $2$-BSEC$(3,c,3,2)$ for any $c\geq 3$ and $c\neq 4$.
\end{Lemma}

\proof When $c\in\{3,5,7,9,11,13,15,17,19,25\}$, take two copies of a $2$-BSEC$(3,c,3,1)$ from Theorem \ref{3-BSEC}(1). When $c\in\{12,18\}$, the conclusion follows from Theorem \ref{3-BSEC}(2). When $c\in\{6,8\}$, see Appendix C. When $c\in\{10,14,16,20,22\}$ see Appendix E. When $c\geq 21$ and $c\neq 22,25$, use induction on $c$ and apply Lemma \ref{3c+x recur} with $x=3,5,7$. \qed

\begin{Lemma}\label{m=6 index 2}
There exists a $2$-BSEC$(6,c,3,2)$ for any $c\geq 3$.
\end{Lemma}

\proof When $c=3$, the conclusion follows from Lemma \ref{m=3 index 2}. When $c\in\{9,12,15,18\}$, the conclusion follows from Theorem \ref{3-BSEC}(2). When $c\in\{4,6,7,8,10,13,16,19\}$, see Appendix C. When $c\in\{11,14,17,20\}$, see Appendix D. When $c=5$, see Appendix E. When $c\geq 21$, use induction on $c$ and apply Lemma \ref{3c+x recur} with $x=3,4,5$. \qed

\begin{Lemma}\label{index two-3}
Let $r\equiv 0\ ({\rm mod}\ 3)$. There exists a $2$-BSEC$(r,c,3,2)$ for any $c\geq 3$  and $c\neq 4$.
\end{Lemma}

\proof When $r=3,6$, the conclusion follows by Lemmas \ref{m=3 index 2} and \ref{m=6 index 2}. When $r\equiv 0\ ({\rm mod}\ 3)$, $r \geq9$ and $c\geq 3$, $c\neq 4$, by Theorem  \ref{3-HGDD-wei} there is a $(3,2)$-HGDD of type $(c,3^{r/3})$. Apply Construction \ref{construction from hgdd} to obtain a $2$-BSEC$(c,r,3,2)$ (i.e., a $2$-BSEC$(r,c,3,2)$), where the needed $(3,2)$-QMGDD of type $3^c$ (i.e., a $2$-BSEC$(c,3,3,2)$) is from Lemma \ref{m=3 index 2}, and the needed 1-BSEC$(r,3,2)$ comes from Theorem \ref{one-dimensional BSEC}. \qed

\begin{Lemma}\label{index two-3-1}
Let $r\equiv 1\ ({\rm mod}\ 3)$ and $r>4$. If there exists a $2$-BSEC$(r,s,3,2)$ for $s\in\{5,8,11\}$, then there exists a $2$-BSEC$(r,c,3,2)$ for any $c\equiv 2\ ({\rm mod}\ 3)$ and $c\geq 41$.
\end{Lemma}

\proof Let $a\geq 4$ and $c=9a+s$, $s\in\{5,8,11\}$. Lemma \ref{index two-3} implies that a 2-BSEC$(r,3a,3,2)$ exists for any $r\equiv 1\ ({\rm mod}\ 3)$ and $r>4$. Then apply Lemma \ref{3c+x recur} with the given 2-BSEC$(r,s,3,2)$ to obtain a 2-BSEC$(r,9a+s,3,2)$. \qed

\begin{Lemma}\label{index two-3-2}
Let $c\equiv 2\ ({\rm mod}\ 3)$ and $5\leq c\leq 38$. If there exists a $2$-BSEC$(s,c,3,2)$ for $s\in\{4,7,10\}$, then there exists a $2$-BSEC$(r,c,3,2)$ for any $r\equiv 1\ ({\rm mod}\ 3)$ and $r\geq 40$.
\end{Lemma}

\proof Let $a\geq 4$ and $r=9a+s$, $s\in\{4,7,10\}$. Lemma \ref{index two-3} shows that a 2-BSEC$(3a,c,3,2)$ exists for any $c\equiv 2\ ({\rm mod}\ 3)$ and $5\leq c\leq 38$. Then apply Lemma \ref{3c+x recur} with the given 2-BSEC$(s,c,3,2)$ to obtain a 2-BSEC$(9a+s,c,3,2)$. \qed

Combining the results of Lemmas \ref{index two-3-1} and \ref{index two-3-2}, we can see that if there exists a 2-BSEC$(r,c,3,2)$ for $r\equiv 1\ ({\rm mod}\ 3)$, $4\leq r\leq 37$ and $c\equiv 2\ ({\rm mod}\ 3)$, $5\leq c\leq 38$, then there exists a 2-BSEC$(r,c,3,2)$ for any $r\equiv 1\ ({\rm mod}\ 3)$, $r> 4$ and $c\equiv 2\ ({\rm mod}\ 3)$, $c\geq 5$.
Thus we would almost complete the proof for the existence of 2-BSEC$(r,c,3,2)$s if we could find all small orders for admissible $r,c\geq 4$ and $r,c\leq 38$.

Even if by computer search we can find all possible examples for these small orders, it is not a good way to write them down since they would occupy too many pages, and especially, it would be uninteresting and ugly. Here we will provide another proof that depends on direct constructions via sequences. The new proof can provide 2-BSECs admitting good algebraic structures and have clear advantage over those with no algebraic structures in the identification of the supports.

\begin{Lemma}{\rm (Bryant et al. 2002)}\label{sequence}
There exist triples $T_1,T_2,\ldots,T_x$ that partition:
\begin{enumerate}
\item [$(1)$] either $\{3,4,\ldots,3x+2\}$ or $\{3,4,\ldots,3x+1,3x+3\}$ if $v=6x+5$ and $x\geq 2$;
\item [$(2)$] either $\{4,5,\ldots,3x+3\}$ or $\{4,5,\ldots,3x+2,3x+4\}$ if $v=6x+7$ and $x\geq 3$,
\end{enumerate}
so that if $T_i=\{a_i,b_i,c_i\}$ then either $a_i+b_i=c_i$ or $a_i+b_i+c_i=v$.
\end{Lemma}

In this paper, we always assume that $[a,b]$ denotes the set of integers $n$ such that $a\leq n\leq b$.

\begin{Lemma}
\label{se}{\rm (Zhang and Chang 2005)} For $1\leq d\leq4,$ if
$(m,k)\equiv (0,1),(1,d),(2,0),(3,d+1)\ ({\rm mod}\ (4,2))$ such
that $m\geq 2d-3$ and $(m/2)(2d-1-m)+1\leq k\leq (m/2)(m-2d+5)+1,$
then $[d,d+3m]\backslash\{k+d+m-1\}$ can be partitioned into triples
$\{a_i,b_i,c_i\}, 1\leq i\leq m,$ such that $a_i+b_i=c_i.$
\end{Lemma}

\begin{Lemma}\label{m=4 index 2}
There exists a $2$-BSEC$(r,4,3,2)$ for any $r\equiv 2\ ({\rm mod}\ 3)$ and $r\geq 5$.
\end{Lemma}

\proof For $r\equiv 5\ ({\rm mod}\ 6)$, let $r=6x+5$ and $x \geq0$. When $x=1$, a $2$-BSEC$(11,4,3,2)$ follows from Example \ref{11,4,3,2}.
When $x\geq 2$, by Lemma \ref{sequence}(1), there exists a partition of $\{3,4,\ldots,3x+2\}$ or $\{3,4,\ldots,3x+1,3x+3\}$ into triples $T_1,T_2,\ldots,T_{x}$ such that for $1\leq i\leq x$, the elements in $T_i=\{a_i,b_i,c_i\}$ are named so that $a_i+b_i=c_i$ or $a_i+b_i+c_i=6x+5$. For $1\leq i\leq x$ and $j\in Z_r$, let $T_i(j)=\{j,(3x+2)a_i+j,(3x+2)(a_i+b_i)+j\}$, reducing the sums modulo $6x+5$. Let ${\cal T}=\{T_i(j):1\leq i\leq x,j\in Z_r\}$. For each triple $T\in{\cal T}$, construct a $(3,1)$-GDD of type $4^3$ on the set $T\times Z_4$ with group set $\{\{l\}\times Z_4:l\in T\}$; this GDD exists by Theorem \ref{3-GDD}. Denote by ${\cal A}_T$ the set of its blocks. Let ${\cal B}_1=\bigcup_{T\in {\cal T}}{\cal A}_T$.

When $x=0$ or $x\geq2$, ${\cal B}_2$ consists of $20r$ blocks which can be obtained from the following 5 base blocks by $(+1\ {\rm mod}\ r,+1\ {\rm mod}\ 4)$:
\begin{center}
\begin{tabular}{lll}
 $\{(0,0),(0,2),(6x+4,3)\}$, & $\{(0,0),(3x+2,0),(6x+4,2)\}$, \\
 $\{(0,0),(3x+2,1),(6x+4,2)\}$, & $\{(0,0),(3x+2,0),(6x+4,3)\}$,\\
 $\{(0,0),(3x+2,2),(6x+4,1)\}$. &
\end{tabular}
\end{center}
Denote by $\cal D$ the set of the above 5 base blocks. Then $\bigcup_{D\in{\cal D}}\Delta(D)=2\cdot\{(0,2),(\pm 1,\pm 1),\linebreak(\pm 1,2),(\pm(3x+2),0),(\pm (3x+2),\pm 1),(\pm(3x+2),2)\}$. It is readily checked that $(2\cdot{\cal B}_1)\cup{\cal B}_2$ forms the required $2$-BSEC$(r,4,3,2)$.

For $r\equiv 2\ ({\rm mod}\ 6)$ and $r\geq 8$, ${\cal B}_3$ consists of $4r^2$ blocks which can be obtained from the following $r$ base blocks by $(+1\ {\rm mod}\ r,+1\ {\rm mod}\ 4)$:
\begin{center}
\begin{tabular}{lll}
 $\{(0,0),(i,1),(2i+1,2)\}$, & $i\in[1,r/2]$; \\
 $\{(0,0),(i,1),(2i,2)\}$, & $i\in[r/2+2,r-1]$;
\end{tabular}
\begin{tabular}{lll}
 $\{(0,0),(r/2,0),(1,1)\}$, & $\{(0,0),(2,0),(2,2)\}$.
\end{tabular}
\end{center}
By Lemma \ref{se}, when $r\equiv 2\ ({\rm mod}\ 12)$ and $r>2$, $[2,r-3]\setminus \{r/2\}$ (when $r\equiv 8\ ({\rm mod}\ 12)$ and $r>8$, $[3,r-2]\setminus \{r/2\}$) can be partitioned into triples $\{a_i,b_i,c_i\}$, such that $a_i+b_i=c_i$, $i\in[1,(r-5)/3]$. When $r=8$, $[3,r-2]\setminus \{r/2\}=\{3,5,6\}$ can be partitioned into one triple $\{a_1,b_1,c_1\}=\{5,6,3\}$, such that $a_1+b_1\equiv c_1\ ({\rm mod} \ 8)$. Let ${\cal B}_4$ consist of $4r(r-5)/3$ blocks which can be obtained from the following $(r-5)/3$ base blocks by $(+1\ {\rm mod}\ r,+1\ {\rm mod}\ 4)$:
$$\{(0,0),(a_i,0),(c_i,0)\},\ \ i\in[1,(r-5)/3].$$
Then ${\cal B}_3\cup{\cal B}_4$ forms the required $2$-BSEC$(r,4,3,2)$. \qed

\begin{Lemma}\label{index two-n1}
There exists a $2$-BSEC$(r,c,3,2)$ for any $r\equiv 2\ ({\rm mod}\ 3)$, $r\geq 5$ and $c\equiv 4\ ({\rm mod}\ 12)$.
\end{Lemma}

\proof When $c=4$, the conclusion follows from Lemma \ref{m=4 index 2}. When $c\equiv 4\ ({\rm mod}\ 12)$ and $c\geq 16$, take a $(3,2)$-IGDD of type $(4,2)^{c/4}$ from Theorem \ref{3-GDD}. By Theorem \ref{MGDD}(1), there exists a (3,1)-MGDD of type $3^r$ for any $r\equiv 2\ ({\rm mod}\ 3)$ and $r\geq 5$. Then apply Construction \ref{Construction from IGDD} to obtain a 2-BSEC$(r,c,3,2)$, where the needed $(3,2)$-GDD of type $2^{c/4}$ is from Theorem \ref{3-GDD}, and the needed 2-BSEC$(r,4,3,2)$ is from Lemma \ref{m=4 index 2}. \qed

\begin{Lemma}\label{7 for lambda2}
There exists a $2$-BSEC$(7,c,3,2)$ for any $c\equiv 2\ ({\rm mod}\ 6)$ and $c\geq 8$.
\end{Lemma}

\proof When $c=8$, a $2$-BSEC$(7,8,3,2)$ can be found in Appendix A.
Assume that $c\equiv 2\ ({\rm mod}\ 6)$ and $c\geq 14$. Let ${\cal T}=\{\{j,1+j,3+j\}:j\in Z_7\}$, reducing the sums modulo 7. For each triple $T\in{\cal T}$, construct a $(3,1)$-MGDD of type $c^3$ on the set $T\times Z_c$ with group set $\{\{l\}\times Z_c:l\in T\}$ and hole set $\{T\times\{l\}:l\in Z_c\}$; this MGDD exists by Theorem \ref{MGDD}(1). Denote by ${\cal A}_T$ the set of its blocks. Let ${\cal B}_1=\bigcup_{T\in {\cal T}}{\cal A}_T$.

${\cal B}_2$ consists of $7c(c+1)$ blocks which can be obtained from the following $c+1$ base blocks by $(+1\ {\rm mod}\ 7,+1\ {\rm mod}\ c)$.
\begin{center}
\begin{tabular}{lll}
$\{(0,0),(1,i),(3,2i+1)\}$, & $i\in[1,c/2]$;  \\
$\{(0,0),(1,i),(3,2i+2)\}$, & $i\in[c/2+1,c-1]$;
\end{tabular}
\begin{tabular}{lll}
 $\{(0,0),(0,c/2+2),(2,c/2+2)\}$, &  $\{(0,0),(0,2),(3,2)\}$.
\end{tabular}
\end{center}
By Lemma \ref{se}, when $c\equiv 8\ ({\rm mod}\ 12)$ and $c>8$, $[2,c-3]\setminus \{c/2+2\}$ (when $c\equiv 2\ ({\rm mod}\ 12)$ and $c>14$, $[3,c-2]\setminus \{c/2+2\}$) can be partitioned into triples $\{a_i,b_i,c_i\}$, such that $a_i+b_i=c_i$, $i\in[1,(c-5)/3]$. When $c=14$, $[3,c-2]\setminus \{c/2+2\}=[3,12]\setminus \{9\}$ can be partitioned into the triples $\{a_i,b_i,c_i\}$, such that $a_i+b_i\equiv c_i\ ({\rm mod} \ 14)$, $i\in[1,3]$, which are $\{3,7,10\}$, $\{6,12,4\}$ and $\{8,11,5\}$. Let ${\cal B}_3$ consist of $7c(c-5)/3$ blocks which can be obtained from the following $(c-5)/3$ base blocks by $(+1\ {\rm mod}\ 7,+1\ {\rm mod}\ c)$:
$$\{(0,0),(0,a_i),(0,c_i)\},\ \ i\in[1,(c-5)/3].$$
Then ${\cal B}_1\cup{\cal B}_2\cup{\cal B}_3$ forms the required $2$-BSEC$(7,c,3,2)$. \qed

\begin{Lemma}
\label{seq}{\rm (Zhang and Chang 2005)} Let $(m,d)\equiv
(0,1),(1,1),(0,0),(3,0)\ ({\rm mod}\ (4,2))$ such that $m\geq 2d-1.$
Then $[d,d+3m-1]$ can be partitioned into triples $\{a_i,b_i,c_i\},
1\leq i\leq m,$ such that $a_i+b_i=c_i.$
\end{Lemma}

\begin{Lemma}\label{1 mod 6 index 2}
Let $r\equiv 1\ ({\rm mod}\ 6)$, $r\geq 7$, and $c\equiv 2\ ({\rm mod}\ 6)$, $c\geq 8$. There exists a $2$-BSEC$(r,c,3,2)$.
\end{Lemma}

\proof When $r=7$, the conclusion follows from Lemma \ref{7 for lambda2}. Assume that $r\equiv 1\ ({\rm mod}\ 6)$ and $r\geq 13$. By Lemma \ref{seq}, when $r\equiv 1\ ({\rm mod}\ 12)$ and $r\geq13$, $[2,r-3]$ (when $r\equiv 7\ ({\rm mod}\ 12)$ and $r\geq19$, $[3,r-2]$) can be partitioned into triples $\{a_i,b_i,c_i\}$, such that $a_i+b_i=c_i$, $i\in[1,(r-4)/3]$. For $1\leq i\leq (r-4)/3$ and $j\in Z_r$, let $T_i(j)=\{j,a_i+j,c_i+j\}$, reducing the sums modulo $r$. Let ${\cal T}=\{T_i(j):1\leq i\leq (r-4)/3,j\in Z_r\}$. For each triple $T\in{\cal T}$, construct a $(3,1)$-GDD of type $c^3$ on the set $T\times Z_c$ with group set $\{\{l\}\times Z_c:l\in T\}$; this GDD exists by Theorem \ref{3-GDD}. Denote by ${\cal A}_T$ the set of its blocks. Let ${\cal B}_1=\bigcup_{T\in {\cal T}}{\cal A}_T$.

${\cal B}_2$ consists of $rc^2$ blocks which can be obtained from the following $c$ base blocks by $(+1\ {\rm mod}\ r,+1\ {\rm mod}\ c)$.
\begin{center}
\begin{tabular}{lll}
$\{(0,0),(1,i),(2,1+2i)\}$, & $i\in[1,c/2]$;  \\
$\{(0,0),(1,i),(2,2i)\}$, & $i\in[c/2+2,c-1]$;
\end{tabular}
\begin{tabular}{lll}
$\{(0,0),(0,c/2),(1,1)\}$, & $\{(0,0),(0,c-2),(2,0)\}$.
\end{tabular}
\end{center}
By Lemma \ref{se}, when $c\equiv 8\ ({\rm mod}\ 12)$, $[2,c-3]\setminus \{c/2\}$ (when $c\equiv 2\ ({\rm mod}\ 12)$ and $c>14$, $[3,c-2]\setminus \{c/2\}$) can be partitioned into triples $\{a_i,b_i,c_i\}$, such that $a_i+b_i=c_i$, $i\in[1,(c-5)/3]$. For $c=14$, $[3,c-2]\setminus \{c/2\}=[3,12]\setminus \{7\}$ can be partitioned into the triples $\{a_i,b_i,c_i\}$, such that $a_i+b_i\equiv c_i\ ({\rm mod} \ 14)$, $i\in[1,3]$, which are $\{6,12,4\}$, $\{9,10,5\}$, $\{3,8,11\}$. Let ${\cal B}_3$ consist of $rc(c-5)/3$ blocks which can be obtained from the following $(c-5)/3$ base blocks by $(+1\ {\rm mod}\ r,+1\ {\rm mod}\ c)$:
$$\{(0,0),(0,a_i),(0,c_i)\},\ \ i\in[1,(c-5)/3].$$
Then ${\cal B}_1\cup{\cal B}_2\cup{\cal B}_3$ forms the required $2$-BSEC$(r,c,3,2)$. \qed

\begin{Lemma}\label{10 mod 12 index 2}
There exists a $2$-BSEC$(r,10,3,2)$ for any $r\equiv 2\ ({\rm mod}\ 3)$ and $r\geq 5$.
\end{Lemma}

\proof  When $r=5,8,26$, see Appendix E, A and C, respectively. When $r\in\{11,14,17,20\}$, see Appendix D. When $r=23$, or $r\equiv 2\ ({\rm mod}\ 3)$ and $r\geq 29$, use induction on $r$ and apply Lemmas \ref{3c+x recur} and \ref{index two-3} with $x=5,8$. \qed

\begin{Lemma}\label{4,c,3,2}
There exists a $2$-BSEC$(4,c,3,2)$ for any $c\equiv 0\ ({\rm mod}\ 3)$ and $c\geq 6$.
\end{Lemma}

\proof When $c=6$, the conclusion follows from Lemma \ref{m=6 index 2}. When $c\in\{9,12,15,18,21,27\}$, see Appendix B. When $c=24$, or $c\equiv 0\ ({\rm mod}\ 3)$ and $c\geq30$, use induction on $c$ and apply Lemmas \ref{3c+x recur} and \ref{m=4 index 2} with $x=6,9$. \qed

\begin{Lemma}\label{index two}
There exists a $2$-BSEC$(r,c,3,2)$ for any $rc\equiv 0,2\ ({\rm mod}\ 3)$ and $r,c\geq3$ except when $(r,c)\in\{(3,4),(4,3)\}$.
\end{Lemma}

\proof When $r\equiv 0\ ({\rm mod}\ 3)$, the conclusion follows from Lemmas \ref{nonexistence}, \ref{index two-3} and \ref{4,c,3,2}. When $r\equiv 1\ ({\rm mod}\ 6)$ and $c\equiv 2\ ({\rm mod}\ 6)$, the conclusion follows from Lemma \ref{1 mod 6 index 2}. When $r\equiv 1\ ({\rm mod}\ 6)$ and $c\equiv 5\ ({\rm mod}\ 6)$, take two copies of a $2$-BSEC$(r,c,3,1)$ from Theorem \ref{3-BSEC}(1). When $r\equiv 4\ ({\rm mod}\ 12)$ and $c\equiv 2\ ({\rm mod}\ 3)$, the conclusion follows from Lemma \ref{index two-n1}.

When $r\equiv 10\ ({\rm mod}\ 12)$ and $c\equiv 2\ ({\rm mod}\ 3)$, or equivalently, by the symmetry of $r$ and $c$, we consider the case $r\equiv 2\ ({\rm mod}\ 3)$ and $c\equiv 10\ ({\rm mod}\ 12)$. Let $s\in\{4,7,10\}$. Let $a\equiv s+2\ ({\rm mod}\ 12)$ and $a\geq s+2$. So $a\equiv 6,9$ or $12\ ({\rm mod}\ 12)$. Write $c=3a+s$. Then $c\equiv 22,34$ or $46\ ({\rm mod}\ 36)$ and $c\geq 22,34$ or $46$, respectively. Lemma \ref{index two-3} shows that a 2-BSEC$(r,a,3,2)$ exists for any $r\equiv 2\ ({\rm mod}\ 3)$ and $r\geq 5$. By Lemmas \ref{m=4 index 2}, \ref{7 for lambda2} and \ref{10 mod 12 index 2}, there exists a 2-BSEC$(r,s,3,2)$ for $s\in\{4,7,10\}$. Now apply Lemma \ref{3c+x recur} to obtain a 2-BSEC$(r,3a+s,3,2)$. \qed

\subsection{$\lambda=3$}

\begin{Lemma}\label{m=5 index three}
There exists a $2$-BSEC$(5,c,3,3)$ for any $c\equiv 5\ ({\rm mod}\ 6)$.
\end{Lemma}

\proof When $c=5$, see Appendix A.
When $c\equiv 5\ ({\rm mod}\ 6)$ and $c\geq 11$, start from a $(3,3)$-HGDD of type $(5,3^{(c-2)/3} 2^1)$, which is from Theorem \ref{3-HGDD} by repeating blocks. Apply Construction \ref{construction from hgdd} to obtain a $2$-BSEC$(5,c,3,3)$, where the needed $(3,3)$-QMGDD of type $2^5$ is from Example \ref{example QMGDD-(2^5)} by repeating blocks, the needed $(3,3)$-QMGDD of type $3^5$ (i.e., a $2$-BSEC$(5,3,3,3)$) is obtained by taking three copies of a $2$-BSEC$(5,3,3,1)$ (from Theorem \ref{3-BSEC}(1)), and the needed 1-BSEC$(c,3,3)$ comes from Theorem \ref{one-dimensional BSEC}. \qed

\begin{Lemma}\label{m=7 index three}
There exists a $2$-BSEC$(7,c,3,3)$ for any $c\equiv 1\ ({\rm mod}\ 6)$ and $c\geq 7$.
\end{Lemma}

\proof When $c=7,13$, see Appendices A and C.
When $c\equiv 1\ ({\rm mod}\ 6)$ and $c\geq 19$, start from a $(3,3)$-HGDD of type $(7,3^{(c-7)/3} 7^1)$, which is from Theorem \ref{3-HGDD} by repeating blocks. Apply Construction \ref{construction from hgdd} to obtain a $2$-BSEC$(7,c,3,3)$, where the needed $(3,3)$-QMGDD of type $7^7$ is from Example \ref{example QMGDD-(2^5)}, the needed $(3,3)$-QMGDD of type $3^7$ (i.e., a $2$-BSEC$(7,3,3,3)$) is obtained by taking three copies of a $2$-BSEC$(7,3,3,1)$ (from Theorem \ref{3-BSEC}(1)),  and the needed 1-BSEC$(c,3,3)$ is from Theorem \ref{one-dimensional BSEC}. \qed

\begin{Lemma}\label{index three}
There exists a $2$-BSEC$(r,c,3,3)$ for any $rc\equiv 1\ ({\rm mod}\ 2)$ and $r,c\geq 3$.
\end{Lemma}

\proof For $r=5$ and $c\equiv 1,3\ ({\rm mod}\ 6)$, or $r=7$ and $c\equiv 3,5\ ({\rm mod}\ 6)$, repeat the blocks of a $2$-BSEC$(r,c,3,1)$ three times (from Theorem \ref{3-BSEC}(1)). For $r=5$ and $c\equiv 5\ ({\rm mod}\ 6)$, or $r=7$ and $c\equiv 1\ ({\rm mod}\ 6)$, the conclusion follows from Lemmas \ref{m=5 index three} and \ref{m=7 index three}, respectively. For the other cases, see Theorem \ref{3-BSEC}(2). \qed

\subsection{$\lambda=6$}

\begin{Lemma}\label{m=4 index six}
There exists a $2$-BSEC$(4,c,3,6)$ for any $c\equiv 1\ ({\rm mod}\ 3)$ and $c\geq 4$.
\end{Lemma}

\proof  When $c=4$, see Appendix A. When $c\equiv 4\ ({\rm mod}\ 12)$ and $c\geq 16$,
by Lemma \ref{se}, $[2,c-2]\setminus\{c/2+2\}$ can be partitioned into triples $\{a_i,b_i,c_i\}$, such that $a_i+b_i=c_i$, $i\in[1,(c-4)/3]$. By Lemma \ref{seq}, $[2,c-3]$ can be partitioned into triples $\{a'_i,b'_i,c'_i\}$, such that $a'_i+b'_i=c'_i$, $i\in[1,(c-4)/3]$. Thus $2\cdot[2,c-3]$ can be partitioned into triples $\{a_i,b_i,c_i\}$, such that $a_i+b_i=c_i$, $i\in[(c-1)/3,c-4]$. Let ${\cal B}_1=\{\{(j,l),(j,a_i+l),(j,c_i+l)\}: 1\leq i\leq c-4, j\in Z_4,l\in Z_c\}$, reducing the sums modulo $c$.

${\cal B}_2$ consists of $4c(3c-1)$ blocks which can be obtained from the following $3c-1$ base blocks by $(+1\ {\rm mod}\ 4,+1\ {\rm mod}\ c)$.
\begin{center}
\begin{tabular}{lll}
$\{(0,0),(1,1+i),(2,2+2i)\}$, & $i\in3\cdot[0,c/2-1]\setminus\{c/2-2,c/2-1\}$;  \\
$\{(0,0),(1,1+i),(2,1+2i)\}$, & $i\in3\cdot[c/2+1,c-2]$;
\end{tabular}
\begin{tabular}{lll}
$\{(0,0),(0,c/2+2),(1,c/2+1)\}$, & $2\{(0,0),(0,c-2),(2,c-1)\}$, \\ $\{(0,0),(1,c/2-1),(2,c-1)\}$,&
$\{(0,0),(1,c/2),(2,1)\}$,  \\ $\{(0,0),(1,c/2-1),(2,0)\}$, & $\{(0,0),(1,c-1),(2,c-2)\}$.
\end{tabular}
\end{center}
Then ${\cal B}_1\cup{\cal B}_2$ forms the required $2$-BSEC$(4,c,3,6)$ for $c\equiv 4\ ({\rm mod}\ 12)$ and $c\geq 16$.

When $c\equiv 10\ ({\rm mod}\ 12)$, by Lemma \ref{se}, $[2,c-2]\setminus\{c/2+2\}$ can be partitioned into triples $\{a_i,b_i,c_i\}$, such that $a_i+b_i=c_i$, $i\in[1,(c-4)/3]$; $[2,c-2]\setminus\{c-3\}$ can be partitioned into triples $\{a'_i,b'_i,c'_i\}$, such that $a'_i+b'_i=c'_i$, $i\in[1,(c-4)/3]$. Thus $2\cdot[2,c-2]\setminus(2\cdot\{c-3\})$ can be partitioned into triples $\{a_i,b_i,c_i\}$, such that $a_i+b_i=c_i$, $i\in[(c-1)/3,c-4]$. Let ${\cal B}_3=\{\{(j,l),(j,a_i+l),(j,c_i+l)\}: 1\leq i\leq c-4, j\in Z_4,l\in Z_c\}$, reducing the sums modulo $c$.

${\cal B}_4$ consists of $4c(3c-1)$ blocks which can be obtained from the following $3c-1$ base blocks by $(+1\ {\rm mod}\ 4,+1\ {\rm mod}\ c)$.
\begin{center}
\begin{tabular}{lll}
$\{(0,0),(1,1+i),(2,2+2i)\}$, & $i\in3\cdot[0,c/2-1]\setminus\{0,c/2-2,c/2-1,c/2-1\}$;  \\
$\{(0,0),(1,1+i),(2,1+2i)\}$, & $i\in3\cdot[c/2+1,c-2]$;
\end{tabular}
\begin{tabular}{lll}
$\{(0,0),(0,c/2+2),(1,c/2+1)\}$, & $\{(0,0),(0,c-3),(2,c-1)\}$, \\ $\{(0,0),(0,c-3),(2,c-2)\}$,&
$2\{(0,0),(1,c/2),(2,1)\}$,  \\ $2\{(0,0),(1,1),(2,0)\}$, & $2\{(0,0),(1,c/2-1),(2,c-1)\}$.
\end{tabular}
\end{center}
Then ${\cal B}_3\cup{\cal B}_4$ forms the required $2$-BSEC$(4,c,3,6)$ for $c\equiv 10\ ({\rm mod}\ 12)$.

When $c\equiv 1\ ({\rm mod}\ 6)$ and $c\geq 7$, let $c=6x+1$ and $x\geq 1$. When $x=2,3$, see Appendix A. When $x\geq4$, by Lemma \ref{sequence}(2), there is a partition of $\{4,5,\ldots,3x\}$ or $\{4,5,\ldots,3x-1,3x+1\}$ into triples $T_1,T_2,\ldots,T_{x-1}$ such that for $1\leq i\leq x-1$, the elements in $T_i=\{a_i,b_i,c_i\}$ are named so that $a_i+b_i=c_i$ or $a_i+b_i+c_i=6x+1$. For $1\leq i\leq x-1$ and $j\in Z_c$, let $T_i(j)=\{j,3xa_i+j,3x(a_i+b_i)+j\}$, and $T'_i(j)=\{j,6xa_i+j,6x(a_i+b_i)+j\}$, reducing the sum modulo $c$. Let ${\cal T}=\{T_i(j):1\leq i\leq x-1,j\in Z_c\}$ and ${\cal T}'=\{T'_i(j):1\leq i\leq x-1,j\in Z_c\}$. For $T\in{\cal T}\cup{\cal T}'$, construct a $(3,1)$-GDD of type $4^3$ on the set $Z_4 \times T$ with group set $\{Z_4\times \{l\}:l\in T\}$; this GDD exists by Theorem \ref{3-GDD}. Denote by ${\cal A}_T$ the set of its blocks. Let ${\cal B}_5=\bigcup_{T\in {\cal T}}{\cal A}_T$ and ${\cal B}'_5=\bigcup_{T\in {\cal T}'}{\cal A}_T$.

When $x=1$ or $x\geq4$, ${\cal B}_6$ consists of $92c$ blocks which can be obtained from the following 23 base blocks by $(+1\ {\rm mod}\ 4,+1\ {\rm mod}\ c)$:
\begin{center}\tabcolsep 0.05in
\begin{tabular}{lll}
$\{(0,0),(2,0),(3,3x)\}$, & $\{(0,0),(2,0),(2,3x-1)\}$, & $\{(0,0),(2,0),(3,3x-1)\}$,\\
$\{(0,0),(1,3x),(0,6x-1)\}$, & $\{(0,0),(0,3x-1),(2,6x-2)\}$, & $\{(0,0),(1,3x),(3,6x)\}$,\\
$\{(0,0),(2,3x),(1,6x)\}$, & $\{(0,0),(1,6x),(0,6x-2)\}$, & $\{(0,0),(2,6x),(3,6x-2)\}$,\\
$\{(0,0),(3,6x),(1,6x-2)\}$, & $2\{(0,0),(2,3x),(0,3x-1)\}$, & $\{(0,0),(0,3x),(1,3x-1)\}$,\\
$2\{(0,0),(3,3x),(1,3x-1)\}$, & $2\{(0,0),(1,3x),(2,3x-1)\}$, & $3\{(0,0),(0,3x),(3,3x-1)\}$,\\
$\{(0,0),(3,3x),(0,3x-1)\}$, & $\{(0,0),(0,3x),(2,3x-1)\}$, & $\{(0,0),(2,3x),(1,3x-1)\}$.
\end{tabular}
\end{center}
Then $(5\cdot{\cal B}_5)\cup{\cal B}'_5\cup{\cal B}_6$ forms the required $2$-BSEC$(4,c,3,6)$.  \qed

\begin{Lemma}\label{m=5 index six}
There exists a $2$-BSEC$(5,c,3,6)$ for any $c\equiv 2\ ({\rm mod}\ 6)$ and $c\geq 8$.
\end{Lemma}

\proof When $c=8$, see Appendix A. When $c\equiv 2\ ({\rm mod}\ 6)$ and $c\geq 14$, take a $(3,6)$-HGDD of type $(5,2^{c/2})$ from Theorem \ref{3-HGDD-wei}. Then apply Construction \ref{construction from hgdd} to obtain a $2$-BSEC$(5,c,3,6)$, where the needed $(3,6)$-QMGDD of type $2^5$ is from Example \ref{example QMGDD-(2^5)} by repeating blocks, and the needed $1$-BSEC$(c,3,6)$ is from Theorem \ref{one-dimensional BSEC}. \qed

\begin{Lemma}\label{m=7 index six}
There exists a $2$-BSEC$(7,c,3,6)$ for any $c\equiv 4\ ({\rm mod}\ 6)$.
\end{Lemma}

\proof When $c=4$, the conclusion follows from Lemma \ref{m=4 index six}. When $c\equiv 4\ ({\rm mod}\ 6)$ and $c\geq 10$, take a $(3,6)$-HGDD of type $(7,2^{c/2})$ from Theorem \ref{3-HGDD-wei}. Then apply Construction \ref{construction from hgdd} to obtain a $2$-BSEC$(7,c,3,6)$, where the needed $(3,6)$-QMGDD of type $2^7$ is from Example \ref{example QMGDD-(2^5)}, and the needed $1$-BSEC$(c,3,6)$ is from Theorem \ref{one-dimensional BSEC}. \qed

\begin{Lemma}\label{m=8 index six}
There exists a $2$-BSEC$(8,c,3,6)$ for any $c\equiv 2\ ({\rm mod}\ 3)$ and $c\geq 5$.
\end{Lemma}

\proof When $c=5,8,11$, see Appendices A and C. When $c\equiv 2\ ({\rm mod}\ 6)$ and $c\geq 14$, take a $(3,6)$-HGDD of type $(8,2^{c/2})$ from Theorem \ref{3-HGDD-wei}. Then apply Construction \ref{construction from hgdd} to obtain a $2$-BSEC$(8,c,3,6)$, where the needed $(3,6)$-QMGDD of type $2^8$ is from Lemma \ref{small QMGDD} by repeating blocks, and the needed $1$-BSEC$(c,3,6)$ is from Theorem \ref{one-dimensional BSEC}.

When $c\equiv 5\ ({\rm mod}\ 6)$ and $c\geq 17$, take a $(3,6)$-HGDD of type $(8,3^{(c-5)/3} 5^1)$ from Theorem \ref{3-HGDD}. Apply Construction \ref{construction from hgdd} to obtain a $2$-BSEC$(8,c,3,6)$, where the needed $(3,6)$-QMGDD of type $5^8$ is from Example \ref{example QMGDD-(2^5)} by repeating blocks, the needed $(3,6)$-QMGDD of type $3^8$ (i.e., a $2$-BSEC$(8,3,3,6)$) is obtained by taking three copies of a $2$-BSEC$(8,3,3,2)$ (from Lemma \ref{index two}), and the needed 1-BSEC$(c,3,6)$ comes from Theorem \ref{one-dimensional BSEC}. \qed

\begin{Lemma}\label{index six}
There exists a $2$-BSEC$(r,c,3,6)$ for any $r,c\geq3$.
\end{Lemma}

\proof When $r=4$ and $c\equiv 0,2\ ({\rm mod}\ 3)$, or $r=5$ and $c\equiv 0,4\ ({\rm mod}\ 6)$, or $r=7$ and $c\equiv 0,2\ ({\rm mod}\ 6)$, or $r=8$ and $c\equiv 0,1\ ({\rm mod}\ 3)$, or $r=6$, repeat the blocks of a $2$-BSEC$(r,c,3,2)$ three times (from Lemma \ref{index two}). When $r=4$ and $c\equiv 1\ ({\rm mod}\ 3)$, or $r=5$ and $c\equiv 2\ ({\rm mod}\ 6)$, or $r=7$ and $c\equiv 4\ ({\rm mod}\ 6)$, or $r=8$ and $c\equiv 2\ ({\rm mod}\ 3)$, the conclusion follows from Lemmas \ref{m=4 index six}, \ref{m=5 index six}, \ref{m=7 index six} and \ref{m=8 index six}, respectively. When $r=5,7$ and $c\equiv 1\ ({\rm mod}\ 2)$, repeat the blocks of a $2$-BSEC$(r,c,3,3)$ twice (from Lemma \ref{index three}). For the other cases, see Theorem \ref{3-BSEC}(2). \qed

\subsection{General $\lambda$}

\begin{Theorem}\label{main result}
Let $r,c\geq3$. There exists a $2$-BSEC$(r,c,3,\lambda)$ if and only if $\lambda rc(rc-5)\equiv 0\ ({\rm mod}\ 6)$, and $\lambda (rc-5)\equiv 0\ ({\rm mod}\ 2)$ except when $(r,c)\in\{(3,4),(4,3)\}$. No $2$-BSEC$(3,4,3,\lambda)$ exists for any $\lambda$.
\end{Theorem}

\proof The necessity comes from Lemma \ref{necessity for any index}. The nonexistence is from Lemma \ref{nonexistence}. When $\lambda\equiv 1,5\ ({\rm mod}\ 6)$, $rc\equiv 1\ ({\rm mod}\ 2)$ and $rc(rc-5)\equiv 0\ ({\rm mod}\ 6)$, repeat the blocks of a $2$-BSEC$(r,c,3,1)$ $\lambda$ times (from Theorem \ref{3-BSEC}(1)).
When $\lambda\equiv 2,4\ ({\rm mod}\ 6)$,  $rc\equiv 0,2\ ({\rm mod}\ 3)$, repeat the blocks of a $2$-BSEC$(r,c,3,2)$ $\lambda/2$ times (from Lemma \ref{index two}).
When $\lambda\equiv 3\ ({\rm mod}\ 6)$, $rc\equiv 1\ ({\rm mod}\ 2)$, repeat the blocks of a $2$-BSEC$(r,c,3,3)$ $\lambda/3$ times (from Lemma \ref{index three}). When $\lambda\equiv 0\ ({\rm mod}\ 6)$, $r$ and $c$ are positive integers, repeat the blocks of a $2$-BSEC$(r,c,3,6)$ $\lambda/6$ times from Lemma \ref{index six}. \qed

\section{2-BSAs with adjacency scheme ``Island''}

In this section we initial the study on 2-BSAs with adjacency scheme ``Island''. For $(x,y)\in Z_r\times Z_c$, the points $(x-1,y)$, $(x+1,y)$, $(x,y-1)$, $(x,y+1)$, $(x-1,y-1)$, $(x-1,y+1)$, $(x+1,y-1)$ and $(x+1,y+1)$ (reducing the arithmetics modulo $r$ and $c$ in the first and second coordinates, respectively) are said to be {\em strongly $2$-contiguous} to the point $(x, y)$.

A {\em two-dimensional balanced sampling plan avoiding strongly $2$-contiguous units} is a pair $(X,{\cal B})$, where $X=Z_r\times Z_c$ and $\cal B$ is a collection of $k$-subsets of $X$ (called {\em blocks}) such that any two strongly $2$-contiguous points do not appear in any block while any two points that are not strongly $2$-contiguous appear in exactly $\lambda$ blocks. It is denoted by a 2-BSA$(r,c,k,\lambda;{\rm IS})$.

By counting the number of blocks and the number of blocks containing a given point in a 2-BSA$(r,c,k,\lambda;{\rm IS})$, a 2-BSA$(r,c,k,\lambda;{\rm IS})$ exists only if $\lambda rc(rc-9)\equiv 0\ ({\rm mod}\ 6)$, and $\lambda (rc-9)\equiv 0\ ({\rm mod}\ 2)$.

The following result is straightforward.

\begin{Lemma}\label{difference method-Is}
Suppose there exist $k$-subsets $B_1,B_2,\ldots,B_b$ of $Z_r\times Z_c$ such that
 $$\bigcup_{i=1}^b \Delta B_i=\lambda\cdot(Z_r\times Z_c\setminus \{(0,0),(0,\pm 1),(\pm1,0),(\pm1,\pm1)\}).$$
Then there exists a $2$-BSA$(r,c,k,\lambda;{\rm IS})$.
\end{Lemma}

The subsets $B_1,B_2,\ldots,B_b$ of $Z_r\times Z_c$ in Lemma \ref{difference method-Is} are called {\em base blocks} of the 2-BSA.

\begin{Example}\label{3,9,3,1}
There is a $2$-BSA$(3,9,3,1;{\rm IS})$. The three base blocks are
$\{(0,0),(0,2),(1,6)\}$, $\{(0,0),(0,3),(1,5)\}$ and $\{(0,0),(0,4),(1,7)\}$.
\end{Example}

In this paper we only establish the existence of a 2-BSA$(3,c,3,1;{\rm IS})$ for odd integer $c$ and $c\geq 9$. Completing the existence of a 2-BSA$(r,c,3,\lambda;{\rm IS})$ will rely heavily on more careful work.

We shall use different strategy not like that in Lemma \ref{difference method-Is}. We require all elements of a 2-BSA$(3,c,3,1;{\rm IS})$ are defined on $Z_{3c}$, and arranged as follows
\begin{center}
$
         \begin{array}{ccccc}
           0 & 1 & 2 & \cdots & c-1 \\
           c & c+1 & c+2 & \cdots & 2c-1 \\
           2c & 2c+1 & 2c+2 & \cdots & 3c-1 \\
         \end{array}
$\end{center}
It is required that every pair $\{x,y\}$ of $Z_{3c}$ satisfying $y-x\in\pm\{0,1,c-1,c,c+1\}$ do not appear in any block. It is easy to see that it is consistent with the definition of a 2-BSA$(3,c,3,1;{\rm IS})$. Thus if the set $[2,(3c-1)/2]\setminus\{c-1,c,c+1\}$ can be partitioned into the triples $\{a_i,b_i,c_i\}$, $i\in[1,(c-3)/2]$, such that $a_i+b_i=c_i$ or $a_i+b_i+c_i\equiv 0\ ({\rm mod}\ 3c)$, then $\{j,a_i+j,a_i+b_i+j\}$, $i\in[1,(c-3)/2]$, $j\in Z_{3c}$, form all blocks of a 2-BSA$(3,c,3,1;{\rm IS})$ for  odd integer $c$ and $c\geq9$.

\begin{Example}\label{3,9,3,1-1}
We give another example of a $2$-BSA$(3,9,3,1;{\rm IS})$. The set $[2,13]\setminus\{8,9,10\}$ can be partitioned into triples $\{2,4,6\}$, $\{5,7,12\}$ and $\{3,11,13\}$. Then $\{j,2+j,6+j\}$, $\{j,5+j,12+j\}$ and $\{j,3+j,14+j\}$, $j\in Z_{27}$,  form all blocks of a $2$-BSA$(3,9,3,1;{\rm IS})$.
\end{Example}

\begin{Theorem}\label{small value}
There exists a $2$-BSA$(3,c,3,1;{\rm IS})$ for any odd integer $c$ and $c\geq 9$.
\end{Theorem}

\proof It suffices to partition the set $[2,(3c-1)/2]\setminus\{c-1,c,c+1\}$ into triples $\{a_i,b_i,c_i\}$, $i\in[1,(c-3)/2]$, such that $a_i+b_i=c_i$ or $a_i+b_i+c_i\equiv 0\ ({\rm mod}\ 3c)$. $c=9$ comes from Example \ref{3,9,3,1-1}.

$\bullet$ $c=11$:
\begin{tabular}{llllll}
$\{2,4,6\}$, & $\{5,9,14\}$, & $\{7,8,15\}$, & $\{3,13,16\}$.
\end{tabular}

$\bullet$ $c=13$:
\begin{tabular}{llllll}
$\{2,4,6\}$, & $\{3,16,19\}$, & $\{5,10,15\}$, & $\{7,11,18\}$, & $\{8,9,17\}$.
\end{tabular}

$\bullet$ $c=15$:
\begin{tabular}{lllllll}
$\{2,10,12\}$, & $\{3,18,21\}$, & $\{5,6,11\}$, & $\{7,13,20\}$, & $\{8,9,17\}$, & $\{4,19,22\}$.
\end{tabular}

$\bullet$ $c=17$:

\begin{tabular}{lllllll}
$\{2,8,10\}$, & $\{3,9,12\}$, & $\{4,19,23\}$, & $\{6,14,20\}$, & $\{7,15,22\}$, & $\{11,13,24\}$, \\ $\{5,21,25\}$.
\end{tabular}

$\bullet$ $c=19$:

\begin{tabular}{lllllll}
$\{2,10,12\}$, & $\{3,24,27\}$, & $\{4,9,13\}$, & $\{5,23,28\}$, & $\{6,16,22\}$, & $\{7,14,21\}$, \\
$\{8,17,25\}$, & $\{11,15,26\}$.&&&&
\end{tabular}

$\bullet$ $c=21$:

\begin{tabular}{lllllll}
$\{2,13,15\}$, & $\{3,26,29\}$, & $\{4,12,16\}$, & $\{5,25,30\}$, & $\{6,11,17\}$, & $\{7,24,31\}$, \\
$\{8,19,27\}$, & $\{9,14,23\}$, & $\{10,18,28\}$.&&&
\end{tabular}

$\bullet$ $c=25$:

\begin{tabular}{lllllll}
$\{2,6,8\}$, & $\{3,20,23\}$, & $\{4,14,18\}$, & $\{5,30,35\}$, & $\{9,27,36\}$, & $\{10,22,32\}$, \\
$\{11,17,28\}$, & $\{12,21,33\}$, & $\{13,16,29\}$, & $\{15,19,34\}$, & $\{7,31,37\}$.&
\end{tabular}

$\bullet$ $c\equiv 1\ ({\rm mod}\ 8)$ and $c\geq 33$:

\begin{tabular}{lll}
$\{4+2i,(3c-11)/4-i,(3c+5)/4+i\}$, & $i\in[0,(c-21)/4]\setminus\{(c-17)/8\}$,\\
$\{5+2i,(5c-5)/4-i,(5c+15)/4+i\}$, & $i\in[0,(c-17)/4]\setminus\{(c-33)/8\}$,
\end{tabular}

\begin{tabular}{lll}
$\{2,(5c+3)/4,(5c+11)/4\}$, & $\{3,(c-9)/2,(c-3)/2\}$,\\
$\{(c-5)/2,(c-1)/2,c-3\}$, & $\{(c+1)/2,(3c-3)/4,(5c-1)/4\}$,\\
$\{(c+3)/2,(c-13)/4,(3c-7)/4\}$, & $\{(3c+1)/4,(5c-5)/8,(11c-3)/8\}$,\\
$\{c+2,(c-1)/4,(5c+7)/4\}$, & $\{c-2,(7c-7)/8,(9c+23)/8\}$.
\end{tabular}

$\bullet$ $c\equiv 3\ ({\rm mod}\ 8)$ and $c\geq27$:

\begin{tabular}{lll}
$\{6+2i,(3c-13)/4-i,(3c+11)/4+i\}$, & $i\in[0,(c-23)/4]\setminus\{(c-19)/8\}$,\\
$\{5+2i,(5c-7)/4-i,(5c+13)/4+i\}$, & $i\in[0,(c-19)/4]\setminus\{(c-27)/8\}$,
\end{tabular}

\begin{tabular}{lll}
$\{2,(3c-9)/4,(3c-1)/4\}$, & $\{3,(3c-5)/4,(3c+7)/4\}$,\\
$\{4,(c-5)/2,(c+3)/2\}$, & $\{(c-7)/4,(c+5)/4,(c-1)/2\}$,\\
$\{(3c+3)/4,(5c-7)/8,(11c-1)/8\}$, & $\{(c-7)/2,(5c+5)/4,(5c+9)/4\}$,\\
$\{(c+1)/2,(5c-3)/4,(5c+1)/4\}$, & $\{(7c+3)/8,(9c+13)/8,c-2\}$,\\
$\{(c-3)/2,c+2,(3c-1)/2\}$.
\end{tabular}

$\bullet$ $c\equiv 5\ ({\rm mod}\ 8)$ and $c\geq29$:

\begin{tabular}{lll}
$\{4+2i,(3c-11)/4-i,(3c+5)/4+i\}$, & $i\in[0,(c-17)/4]\setminus\{(c-29)/8\}$,\\
$\{7+2i,(5c-13)/4-i,(5c+15)/4+i\}$, & $i\in[0,(c-21)/4]\setminus\{(c-21)/8\}$,
\end{tabular}

\begin{tabular}{lll}
$\{2,3,5\}$, & $\{(c-3)/2,(3c-3)/4,(5c-9)/4\}$,\\
$\{(c-1)/2,(3c+1)/4,(5c-1)/4\}$, & $\{(c+1)/2,(3c-7)/4,(5c-5)/4\}$,\\
$\{(c-13)/4,(5c+11)/4,(3c-1)/2\}$, & $\{(c+7)/4,(7c-19)/8,(9c-5)/8\}$,\\
$\{(c-5)/2,(5c+3)/4,(5c+7)/4\}$, & $\{(5c+7)/8,(11c+9)/8,c-2\}$.
\end{tabular}

$\bullet$ $c\equiv 7\ ({\rm mod}\ 8)$ and $c\geq23$:

\begin{tabular}{lll}
$\{4+2i,(3c-9)/4-i,(3c+7)/4+i\}$, & $i\in[0,(c-19)/4]\setminus\{(c-23)/8\}$,\\
$\{7+2i,(5c-11)/4-i,(5c+17)/4+i\}$, & $i\in[0,(c-19)/4]\setminus\{(c-23)/8\}$,
\end{tabular}

\begin{tabular}{lll}
$\{2,3,5\}$, & $\{(c-1)/2,(3c-5)/4,(5c-7)/4\}$,\\
$\{(c+1)/2,(3c+3)/4,(5c+5)/4\}$, & $\{(c-7)/4,(c+3)/2,(3c-1)/4\}$,\\
$\{(c+5)/4,(7c-9)/8,(9c+1)/8\}$, & $\{(c-7)/2,(5c+1)/4,(5c+13)/4\}$,\\
$\{(c-3)/2,(5c-3)/4,(5c+9)/4\}$, & $\{(5c+5)/8,(11c+11)/8,c-2\}$.\\
\end{tabular}

\vspace{0.5cm}

\section{Concluding remarks}

This paper details the three kinds of two-dimensional balanced sampling plans proposed by Wright (2008). A 2-BSA with adjacency scheme ``Row and Column'' is equivalent to a modified group divisible design. 2-BSAs with adjacency scheme ``Sharing a Border'' can be constructed via holey group divisible designs and incomplete group divisible designs. However, it seems that the known recursive constructions can not work for 2-BSAs with adjacency scheme ``Island'', so exploring effective recursive constructions will be interesting but more challenging.

\subsection*{Acknowledgements}
Research of this paper was carried out while the second author was visiting University of Victoria. He expresses his sincere thanks to China Scholarship Council for financial support and to the Department of Mathematics and Statistics in University of
Victoria for the kind hospitality.

\appendix
\section{Appendix}

Here we give some direct constructions for small orders of 2-BSEC$(r,c,3,\lambda)$s on $Z_r\times Z_c$ via Lemma \ref{difference method}. Only base blocks are listed below. All other blocks are obtained by developing these base blocks by $(+1\ {\rm mod}\ r,+1\ {\rm mod}\ c)$.

\vspace{0.2cm}
\noindent(1) $\lambda=2:$

$(r,c)=(7,8):$

\begin{tabular}{llllll}
$\{0_0,1_i,3_{2i}\}$, & $i\in\{2,3,4\}$;  \\
$\{0_0,1_i,3_{2i+1}\}$, & $i\in\{5,6,7\}$; \\
$\{0_0,1_i,3_{2i-1}\}$, & $i\in\{1,2,3,4\}$;  \\
$\{0_0,1_i,3_{2i-2}\}$, & $i\in\{6,7\}$; \\
$\{0_0,1_1,3_0\}$, & $\{0_0,1_5,3_2\}$, & $\{0_0,0_3,2_1\}$, &
 $\{0_0,0_3,3_1\}$, &  $\{0_0,0_2,0_6\}$.
\end{tabular}

$(r,c)=(8,10):$

{\tabcolsep 0.06in
\begin{tabular}{llllll}
$\{0_0,1_i,3_{2i-1}\}$, & $i\in[1,9]\setminus\{5\}$;  \\
$\{0_0,1_i,3_{2i}\}$, & $i\in[1,9]\setminus\{4\}$; \\
$\{0_0,0_{2i+3},4_{5+i}\}$, & $i\in\{0,1,2\}$;\\
$\{0_0,0_2,3_0\}$, & $\{0_0,1_5,4_4\}$, &
$\{0_0,1_4,4_3\}$, & $\{0_0,2_0,4_9\}$, &
$\{0_0,2_4,4_8\}$, & $\{0_0,0_2,0_6\}$.
\end{tabular}}

\noindent (2) $\lambda=3:$

$(r,c)=(5,5):$

\begin{tabular}{llllll}
$\{0_0,0_2,2_3\}$, & $\{0_0,0_2,4_4\}$, &
$\{0_0,0_2,3_2\}$, & $\{0_0,2_0,4_4\}$,&
$\{0_0,2_1,4_3\}$, & $\{0_0,2_2,4_1\}$,\\
$\{0_0,2_0,4_3\}$, & $\{0_0,4_1,3_4\}$,&
$\{0_0,4_2,3_1\}$, & $\{0_0,4_1,3_3\}$.
\end{tabular}

$(r,c)=(7,7):$

\begin{tabular}{lllll}
$\{0_0,l_3,(2l)_{2}\}$, & $l\in\{2,3,4,5,6\}$; \\
$2\{0_0,(3l)_{3},(6l)_{2}\}$, & $l\in\{1,2,3,5,6\}$;\\
$2\{0_0,3_0,3_2\}$, & $2\{0_0,0_3,5_6\}$; &$\{0_0,2_0,2_2\}$,& $\{0_0,0_3,1_6\}$, &
$\{0_0,2_0,4_0\}$.
\end{tabular}

\noindent (3) $\lambda=6:$

$(r,c)=(4,4):$

\begin{tabular}{llllll}
$\{0_0,0_2,1_3\}$, & $\{0_0,0_2,2_2\}$, & $\{0_0,0_2,2_3\}$,&
$\{0_0,1_1,2_0\}$, & $\{0_0,1_1,2_2\}$, & $2\{0_0,1_1,2_3\}$,\\
$\{0_0,1_2,2_0\}$, & $2\{0_0,1_2,2_1\}$, & $\{0_0,1_3,2_2\}$.
\end{tabular}

$(r,c)=(4,13):$

\begin{tabular}{llllll}
$3\{0_0,1_1,2_{12}\}$, & $3\{0_0,1_2,2_{11}\}$, & $3\{0_0,1_3,2_{10}\}$,&
$3\{0_0,1_4,2_9\}$, & $3\{0_0,1_6,2_3\}$, \\ $3\{0_0,1_8,2_7\}$,&
$3\{0_0,1_1,2_8\}$, & $3\{0_0,1_2,2_6\}$, & $2\{0_0,1_3,2_2\}$,&
$3\{0_0,1_5,2_1\}$, \\ $3\{0_0,1_6,2_4\}$, & $2\{0_0,1_8,2_5\}$,&
$\{0_0,1_3,2_0\}$, & $\{0_0,0_4,1_{12}\}$, & $\{0_0,0_2,2_2\}$,\\
$\{0_0,0_5,2_5\}$, & $2\{0_0,0_2,0_5\}$, & $2\{0_0,0_2,0_6\}$,&
$\{0_0,0_3,0_8\}$, & $\{0_0,0_2,0_7\}$,\\ $3\{0_0,0_3,0_7\}$.
\end{tabular}

$(r,c)=(4,19):$

\begin{tabular}{lllll}
$3\{0_0,1_1,2_{18}\}$, & $3\{0_0,1_2,2_{17}\}$, & $3\{0_0,1_3,2_{16}\}$,&
$3\{0_0,1_4,2_{15}\}$, & $3\{0_0,1_5,2_{14}\}$, \\
$3\{0_0,1_6,2_{13}\}$,& $3\{0_0,1_8,2_3\}$, & $3\{0_0,1_{10},2_7\}$, & $3\{0_0,1_{12},2_{11}\}$,&
$3\{0_0,1_{14},2_{12}\}$, \\ $3\{0_0,1_{12},2_9\}$, & $3\{0_0,1_9,2_1\}$,&
$3\{0_0,1_3,2_{10}\}$, & $3\{0_0,1_{10},2_6\}$, & $3\{0_0,1_1,2_5\}$,\\
$3\{0_0,1_5,2_4\}$, & $2\{0_0,1_8,2_2\}$, & $2\{0_0,1_2,2_8\}$,&
$\{0_0,1_6,2_0\}$, & $\{0_0,0_6,1_8\}$, \\ $\{0_0,0_2,2_2\}$,&
$\{0_0,0_8,2_8\}$, & $3\{0_0,0_2,0_7\}$, & $4\{0_0,0_3,0_9\}$,&
$2\{0_0,0_4,0_8\}$, \\ $\{0_0,0_4,0_9\}$, & $\{0_0,0_2,0_4\}$,&
$\{0_0,0_3,0_{11}\}$, & $\{0_0,0_5,0_{12}\}$, & $\{0_0,0_5,0_{11}\}$,\\
$\{0_0,0_3,0_{10}\}$.
\end{tabular}

$(r,c)=(5,8):$

\begin{tabular}{llllll}
$6\{0_0,1_1,2_4\}$, & $5\{0_0,1_2,2_7\}$, & $5\{0_0,1_6,2_5\}$,&
$\{0_0,1_2,2_6\}$, & $\{0_0,1_4,2_0\}$, & $\{0_0,1_4,2_1\}$,\\
$\{0_0,1_4,2_2\}$, & $\{0_0,0_3,1_7\}$, & $4\{0_0,0_2,2_2\}$,&
$2\{0_0,0_2,2_3\}$, & $\{0_0,0_3,2_3\}$, & $2\{0_0,0_3,2_6\}$,\\
$2\{0_0,0_3,2_1\}$, & $\{0_0,0_4,2_5\}$, & $\{0_0,0_4,2_6\}$,&
$\{0_0,0_4,2_7\}$.
\end{tabular}

$(r,c)=(8,8):$

\begin{tabular}{llllll}
$\{0_0,0_2,0_4\}$, & $3\{0_0,0_2,0_5\}$, & $\{0_0,0_2,3_6\}$,&
$\{0_0,0_4,3_6\}$, & $\{0_0,0_4,3_7\}$, & $\{0_0,2_0,4_4\}$,\\
$\{0_0,2_1,4_2\}$, & $\{0_0,2_2,4_7\}$, & $\{0_0,2_3,4_0\}$,&
$\{0_0,2_4,4_6\}$, & $2\{0_0,2_3,4_1\}$, & $\{0_0,2_0,4_7\}$,\\
$\{0_0,2_3,4_2\}$, & $3\{0_0,1_4,3_0\}$, & $3\{0_0,1_1,3_1\}$,&
$3\{0_0,1_2,3_3\}$, & $3\{0_0,1_3,3_5\}$, & $\{0_0,2_3,4_6\}$, \\
$3\{0_0,1_5,3_2\}$, & $3\{0_0,1_6,3_4\}$, & $3\{0_0,1_7,3_6\}$,&
$\{0_0,1_1,3_3\}$, & $\{0_0,1_2,3_0\}$, & $\{0_0,1_3,3_2\}$,\\
$\{0_0,1_4,3_5\}$, & $\{0_0,1_5,3_1\}$, & $\{0_0,1_6,3_3\}$,&
$\{0_0,1_7,3_7\}$, & $\{0_0,1_1,4_0\}$, & $\{0_0,1_2,4_0\}$,\\
$\{0_0,1_2,4_1\}$, & $\{0_0,1_3,4_2\}$, & $\{0_0,1_3,4_5\}$,&
$\{0_0,1_4,4_4\}$, & $\{0_0,1_4,4_5\}$, & $\{0_0,1_1,4_6\}$,\\
$\{0_0,1_5,4_5\}$, & $\{0_0,1_6,4_7\}$, & $\{0_0,1_5,4_4\}$,&
$\{0_0,1_6,4_3\}$, & $2\{0_0,1_7,4_3\}$.
\end{tabular}

\section{Appendix}

For $c\in\{9,12,15,18,21,27\}$, we here explicitly construct a $2$-BSEC$(4,c,3,2)$ as follows.

Let ${\cal B}_1$ consist of $4c/3$ blocks obtained from the $c/3$ blocks $\{(0,i),(0,c/3+i),(0,2c/3+i)\}$, $0\leq i\leq c/3-1$, by $(+1\ {\rm mod}\ 4,-)$. Let  ${\cal B}_2$ consist of $4c(4c/3-2)$ blocks obtained from the following $4c/3-2$ blocks by $(+1\ {\rm mod}\ 4,+1\ {\rm mod}\ c)$.

$\bullet$ $c=9$:

\begin{tabular}{llllll}
$\{0_0,0_2,0_7\}$, & $\{0_0,0_4,1_7\}$, & $\{0_0,0_3,2_8\}$, & $\{0_0,1_8,2_7\}$, & $\{0_0,1_2,2_0\}$,& $\{0_0,1_4,2_1\}$, \\
$\{0_0,1_5,2_2\}$,& $\{0_0,1_1,2_4\}$,& $\{0_0,1_1,2_6\}$,& $\{0_0,1_2,2_6\}$.
\end{tabular}

$\bullet$ $c=12$:

\begin{tabular}{ll}
$\{0_0,1_{1+i},2_{3+2i}\}$, & $i\in[0,5]$,
\end{tabular}

\begin{tabular}{llllll}
$\{0_0,0_3,0_9\}$, & $\{0_0,0_7,0_5\}$, & $\{0_0,0_2,1_{10}\}$, & $\{0_0,0_8,2_{10}\}$, & $\{0_0,1_1,2_0\}$, & $\{0_0,1_7,2_4\}$, \\
$\{0_0,1_8,2_6\}$, & $\{0_0,1_9,2_8\}$.
\end{tabular}

$\bullet$ $c=15$:

\begin{tabular}{ll}
$\{0_0,1_{1+i},2_{3+2i}\}$, & $i\in[0,12]$,
\end{tabular}

\begin{tabular}{llllll}
$\{0_0,0_3,0_{10}\}$, & $\{0_0,0_8,0_4\}$, & $\{0_0,0_9,0_6\}$, & $\{0_0,0_2,1_1\}$, & $\{0_0,0_{13},2_{14}\}$.
\end{tabular}

$\bullet$ $c=18$:

\begin{tabular}{ll}
$\{0_0,1_{1+i},2_{3+2i}\}$, & $i\in[0,8]$,\\
$\{0_0,1_{1+i},2_{4+2i}\}$, & $i\in[9,14]$,
\end{tabular}

\begin{tabular}{llllll}
$\{0_0,0_2,0_{11}\}$, & $\{0_0,0_3,0_{13}\}$, &$\{0_0,0_4,0_{16}\}$, &$\{0_0,0_7,0_{15}\}$,  &$\{0_0,0_5,1_{16}\}$, \\
$\{0_0,0_{14},2_{16}\}$, & $\{0_0,1_{17},2_0\}$.
\end{tabular}

$\bullet$ $c=21$:

\begin{tabular}{ll}
$\{0_0,1_{1+i},2_{3+2i}\}$, & $i\in[0,18]$,
\end{tabular}

\begin{tabular}{llllll}
$\{0_0,0_{11},0_8\}$, & $\{0_0,0_{12},0_6\}$, &$\{0_0,0_{10},0_5\}$, &$\{0_0,0_3,0_{17}\}$,  &$\{0_0,0_4,0_{13}\}$, \\
$\{0_0,0_2,1_1\}$, & $\{0_0,0_{19},2_{20}\}$.
\end{tabular}

$\bullet$ $c=27$:

\begin{tabular}{ll}
$\{0_0,1_{1+i},2_{3+2i}\}$, & $i\in[0,24]$,
\end{tabular}

\begin{tabular}{llllll}
$\{0_0,0_3,0_{15}\}$, & $\{0_0,0_4,0_{22}\}$, &$\{0_0,0_5,0_{21}\}$, &$\{0_0,0_6,0_{20}\}$, &$\{0_0,0_7,0_{24}\}$,\\
$\{0_0,0_8,0_{19}\}$, & $\{0_0,0_{10},0_{23}\}$, & $\{0_0,0_2,1_1\}$, & $\{0_0,0_{25},2_{26}\}$.
\end{tabular}

\noindent Then ${\cal B}_1\cup{\cal B}_2$ forms the required $2$-BSEC$(4,c,3,2)$.

\section{Appendix}

\hspace*{0.5cm} $\bullet$ $2$-BSEC$(3,6,3,2)$

(1) $2\{0_0,0_2,0_4\}$, $2\{0_1,0_3,0_5\}$, $\{0_0,0_3,1_1\}$, $\{0_0,0_3,1_2\}$,
$\{0_1,0_4,1_5\}$, $\{0_1,0_4,1_3\}$, $\{0_2,0_5,1_4\}$, $\{0_2,0_5,1_0\}$, by $(+1\ {\rm mod}\ 3,-)$.

(2)

\begin{tabular}{llllll}
$\{0_0,1_1,2_3\}$, & $\{0_0,1_2,2_4\}$, & $\{0_0,1_3,2_1\}$, & $\{0_0,1_3,2_4\}$, & $\{0_0,1_4,2_5\}$, & $\{0_0,1_4,2_1\}$, \\
$\{0_0,1_5,2_2\}$, & $\{0_0,1_5,2_3\}$, & $\{0_1,1_0,2_5\}$, & $\{0_1,1_0,2_4\}$, &$\{0_1,1_2,2_3\}$, & $\{0_1,1_2,2_5\}$, \\
$\{0_1,1_3,2_4\}$, & $\{0_1,1_4,2_2\}$, & $\{0_1,1_4,2_0\}$, & $\{0_1,1_5,2_2\}$, &$\{0_2,1_0,2_3\}$, & $\{0_2,1_1,2_5\}$, \\
$\{0_2,1_1,2_4\}$, & $\{0_2,1_3,2_5\}$, & $\{0_2,1_3,2_0\}$, & $\{0_2,1_4,2_1\}$, &$\{0_2,1_5,2_1\}$, & $\{0_2,1_5,2_4\}$, \\
$\{0_3,1_0,2_5\}$, & $\{0_3,1_0,2_4\}$, & $\{0_3,1_1,2_2\}$, & $\{0_3,1_2,2_5\}$, &$\{0_3,1_4,2_0\}$, & $\{0_3,1_4,2_2\}$, \\
$\{0_3,1_5,2_1\}$, & $\{0_3,1_5,2_0\}$, & $\{0_4,1_0,2_2\}$, & $\{0_4,1_0,2_3\}$, &$2\{0_4,1_1,2_0\}$, & $2\{0_4,1_2,2_1\}$, \\
$\{0_4,1_3,2_5\}$, & $\{0_4,1_5,2_3\}$, & $\{0_5,1_0,2_1\}$, & $\{0_5,1_1,2_2\}$, &$\{0_5,1_1,2_4\}$, & $\{0_5,1_2,2_0\}$, \\
$\{0_5,1_2,2_3\}$, & $\{0_5,1_3,2_2\}$, & $\{0_5,1_3,2_0\}$, & $\{0_5,1_4,2_3\}$.
\end{tabular}

$\bullet$ $2$-BSEC$(3,8,3,2)$

(1) $\{0_0,0_2,0_4\}$, $\{0_0,1_2,0_5\}$,  $\{0_0,2_2,0_5\}$, by $(+1\ {\rm mod}\ 3,$ $+1\ {\rm mod}\ 8)$;

(2){\tabcolsep 0.05in
\begin{tabular}{llllll}
$\{0_0,1_1,2_5\}$, & $\{2_0,0_1,1_5\}$, &
$\{1_0,2_1,0_2\}$,& $\{0_0,1_1,2_2\}$, &
$\{2_0,1_1,0_2\}$, \\ $\{0_0,2_1,1_4\}$, &
$\{1_0,0_1,2_4\}$, & $\{2_0,1_1,0_4\}$, &
$\{0_0,1_2,2_4\}$,& $\{0_0,2_1,1_3\}$, by $(-,+1\ {\rm mod}\ 8)$.
\end{tabular}}

$\bullet$ $2$-BSEC$(6,6,3,2)$

(1) $\{0_0,0_2,0_4\}$, $\{0_1,0_3,0_5\}$, by $(+1\ {\rm mod}\ 6,-)$.

(2) $\{0_0,0_2,2_3\}$, $\{0_0,0_3,2_2\}$,
$\{0_0,1_3,2_2\}$, $\{0_0,1_4,2_0\}$,
$\{0_0,1_1,3_0\}$,  $\{0_0,1_2,3_3\}$,
$\{0_0,1_3,3_1\}$,  $\{0_0,1_4,3_4\}$,
$\{0_0,1_5,3_2\}$, $\{0_0,1_1,3_5\}$, by $(+1\ {\rm mod}\ 6,$ $+1\ {\rm mod}\ 6)$.

$\bullet$ $2$-BSEC$(6,8,3,2)$

(1) $\{0_0,2_0,4_0\}$, $\{1_0,3_0,5_0\}$, by $(-,+1\ {\rm mod}\ 8)$.

(2) $\{0_0,2_0,0_3\}$, $\{0_0,3_1,2_4\}$, $\{0_0,2_1,1_4\}$, $\{0_0,0_2,5_4\}$,
$\{0_0,2_2,4_4\}$, $\{0_0,3_0,3_4\}$, by $(+1\ \linebreak{\rm mod}\ 6,+1\ {\rm mod}\ 8)$.

(3) $\{0_0,i_1,(2i-1)_3\}$, $i\in[1,5]\setminus\{3\}$; \
$\{0_0,i_1,2i_3\}$, $i\in[1,5]\setminus\{2\}$, by $(+1\ {\rm mod}\ 6,+1\ {\rm mod}\ 8)$.

$\bullet$ $2$-BSEC$(6,c,3,2)$ for $c\in\{4,13,16,19\}$.

Let $\cal A$ consist of $2c$ blocks obtained by developing the following blocks by $(-,+1\ {\rm mod}\ c)$: $\{0_0,2_0,4_0\}$ and $\{1_0,3_0,5_0\}$. Let ${\cal B}$ consist of $36c$ blocks obtained from the following blocks by $(+1\ {\rm mod}\ 6,+1\ {\rm mod}\ c)$

\begin{tabular}{llllll}
$\{0_0,2_1,5_2\}$, & $\{0_0,1_1,4_2\}$, &
$\{0_0,4_1,2_2\}$, & $\{0_0,2_1,1_2\}$, &
$\{0_0,4_0,5_1\}$, & $\{0_0,3_0,3_2\}$.
\end{tabular}

\noindent By Lemma \ref{seq}, when $c=13$, $[2,10]$ (when $c=19$, $[3,17]$) can be partitioned into triples $\{a_i,b_i,c_i\}$, such that $a_i+b_i=c_i$, $i\in[1,(c-4)/3]$. When $c=16$, $[3,14]$ can be partitioned into the triples $\{a_i,b_i,c_i\}$, such that $a_i+b_i\equiv c_i\ ({\rm mod} \ 16)$, $i\in[1,4]$, which are $\{3,5,8\}$, $\{4,7,11\}$, $\{14,12,10\}$ and $\{13,9,6\}$. For $1\leq i\leq (c-4)/3$ and $j\in Z_c$, let $T_i(j)=\{j,a_i+j,c_i+j\}$, reducing the sums modulo $c$. Let ${\cal T}=\{T_i(j):1\leq i\leq (c-4)/3,j\in Z_c\}$. For each triple $T\in{\cal T}$, construct a $(3,1)$-GDD of type $6^3$ on the set $Z_6\times T$ with group set $\{Z_6\times \{l\}:l\in T\}$. Denote by ${\cal A}_T$ the set of its blocks. Let ${\cal C}=\bigcup_{T\in {\cal T}}{\cal A}_T$.
Then ${\cal A}\cup{\cal B}\cup{\cal C}$ forms a $2$-BSEC$(6,c,3,2)$. Note that ${\cal C}=\emptyset$ if $c=4$.

$\bullet$ $2$-BSEC$(6,7,3,2)$.

Let $\cal A$ consist of 14 blocks obtained by developing the following blocks by $(-,+1\ {\rm mod}\ 7)$: $\{0_0,2_0,4_0\}$, $\{1_0,3_0,5_0\}$. Let $\cal B$ consist of 294 blocks obtained by developing the following blocks by $(+1\ {\rm mod}\ 6,+1\ {\rm mod}\ 7)$:

\begin{tabular}{llllll}
$\{0_0,2_0,4_4\}$, & $\{0_0,3_0,4_1\}$, &
$\{0_0,5_1,5_4\}$, & $\{0_0,2_1,5_2\}$, &
$\{0_0,0_2,3_4\}$, & $\{0_0,0_2,1_4\}$, \\
$\{0_0,2_2,0_4\}$.
\end{tabular}

\noindent Let ${\cal T}=\{\{j,1+j,3+j\}:j\in Z_7\}$, reducing the sums modulo $7$. For each triple $T\in{\cal T}$, construct a $(3,1)$-MGDD of type $6^3$ on the set $Z_6\times T$ with group set $\{Z_6\times \{l\}:l\in T\}$ and hole set $\{\{l\}\times T:l\in Z_6\}$. Denote by ${\cal A}_T$ the set of its blocks. Let ${\cal C}=\bigcup_{T\in {\cal T}}{\cal A}_T$. Then ${\cal A}\cup{\cal B}\cup{\cal C}$ forms a $2$-BSEC$(6,7,3,2)$.

$\bullet$ $2$-BSEC$(6,10,3,2)$

Let ${\cal T}=\{\{j,2+j,5+j\}: j\in Z_{10}\}$, reducing the sums modulo $10$. For each triple $T\in{\cal T}$, construct a $(3,1)$-GDD of type $6^3$ on the set $Z_6\times T$ with group set $\{Z_6\times \{l\} :l\in T\}$. Denote by ${\cal A}_T$ the set of its blocks. Let ${\cal B}_1=\bigcup_{T\in {\cal T}}{\cal A}_T$.

Let ${\cal T'}=\{\{j,1+j,4+j\}:j\in Z_{10}\}$, reducing the sums modulo $10$. For each triple $T\in{\cal T'}$, construct a $(3,1)$-MGDD of type $6^3$ on the set $Z_6\times T$ with group set $\{Z_6\times \{l\}:l\in T\}$ and hole set $\{\{l\}\times T:l\in Z_6\}$. Denote by ${\cal A'}_T$ the set of its blocks. Let ${\cal B}_2=\bigcup_{T\in {\cal T'}}{\cal A'}_T$.

${\cal B}_3$ consists of $20$ blocks obtained by developing the following blocks by $(-,+1\ {\rm mod}\ 10)$: $\{0_0,2_0,4_0\}$, $\{1_0,3_0,5_0\}$. ${\cal B}_4$ consists of $420$ blocks obtained by developing the blocks by $(+1\ {\rm mod}\ 6,$ $+1\ {\rm mod}\ 10)$.

\begin{tabular}{llllll}
$\{0_0,2_0,2_4\}$, & $\{0_0,2_1,0_2\}$, &
$\{0_0,1_2,4_4\}$,& $\{0_0,5_1,5_4\}$, &
$\{0_0,3_0,3_4\}$, & $\{0_0,1_1,4_2\}$,\\
$\{0_0,2_2,1_4\}$.
\end{tabular}

\noindent Then $\bigcup_{i=1}^4{\cal B}_i$ forms a $2$-BSEC$(6,10,3,2)$.

$\bullet$ $2$-BSEC$(26,10,3,2)$

Let ${\cal T}=\{\{j,5+j,19+j\},\{j,6+j,17+j\},\{j,7+j,15+j\},\{j,8+j,21+j\},\{j,9+j,12+j\},\{j,10+j,16+j\}: j\in Z_{26}\}$, reducing the sums modulo $26$.  For each triple $T\in{\cal T}$, construct a $(3,1)$-GDD of type $10^3$ on the set $T\times Z_{10}$ with group set $\{\{l\}\times Z_{10}:l\in T\}$. Denote by ${\cal A}_T$ the set of its blocks. Let ${\cal B}_1=\bigcup_{T\in {\cal T}}{\cal A}_T$. ${\cal B}_2$ consists of $25\times 260$ blocks which can be obtained by developing the following $25$ blocks by $(+1\ {\rm mod}\ 26,+1\ {\rm mod}\ 10)$:

\begin{tabular}{llllll}
$\{0_0,1_i,4_{2i-1}\}$, & $i\in[3,8]$;  \\
$\{0_0,2_i,4_{2i}\}$, & $i\in[2,9]$; \\
$\{0_0,0_{2i+3},1_{7+i}\}$, & $i=1,2$;\\
$\{0_0,1_5,2_1\}$, & $\{0_0,1_4,2_1\}$, &
$\{0_0,0_2,4_3\}$, & $\{0_0,1_1,4_9\}$, &
$\{0_0,1_2,4_2\}$, & $\{0_0,0_3,4_0\}$, \\
$\{0_0,1_9,3_9\}$, & $\{0_0,1_1,3_1\}$, &
$\{0_0,0_2,0_6\}$.
\end{tabular}

\noindent Then ${\cal B}_1\cup{\cal B}_2$ forms a $2$-BSEC$(26,10,3,2)$.

$\bullet$ $2$-BSEC$(7,13,3,3)$

Let ${\cal T}=\{\{j,2+j,5+j\},\{j,3+j,7+j\},\{j,2+j,7+j\}: j\in Z_{13}\}$, reducing the sums modulo $13$.  For each triple $T\in{\cal T}$, construct a $(3,1)$-GDD of type $7^3$ on the set $Z_7\times T$ with group set $\{Z_7\times \{l\}:l\in T\}$. Denote by ${\cal A}_T$ the set of its blocks. Let ${\cal B}_1=\bigcup_{T\in {\cal T}}{\cal A}_T$. ${\cal B}_2$ consists of $91\times22$ blocks which can be obtained by developing the following $22$ blocks by $(+1\ {\rm mod}\ 7,+1\ {\rm mod}\ 13)$:

\begin{tabular}{llllll}
$\{0_0,2_0,4_0\}$, & $2\{0_0,3_0,4_1\}$, & $\{0_0,2_0,5_2\}$,&
$\{0_0,2_1,6_3\}$, & $\{0_0,1_1,6_6\}$, & $\{0_0,4_1,0_6\}$,\\
$\{0_0,5_1,2_6\}$, & $\{0_0,0_2,1_6\}$, & $\{0_0,1_2,3_6\}$,&
$\{0_0,2_2,5_6\}$, & $\{0_0,6_2,4_6\}$, & $\{0_0,2_1,6_5\}$,\\
$\{0_0,3_1,2_5\}$, & $\{0_0,6_1,0_5\}$, & $\{0_0,6_1,1_5\}$,&
$\{0_0,2_1,0_4\}$, & $\{0_0,3_1,4_4\}$, & $\{0_0,3_1,6_4\}$,\\
$\{0_0,5_1,5_4\}$, & $\{0_0,5_1,0_4\}$, & $\{0_0,6_1,3_4\}$.
\end{tabular}

\noindent Then ${\cal B}_1\cup{\cal B}_2$ forms a $2$-BSEC$(7,13,3,3)$.

$\bullet$ $2$-BSEC$(8,11,3,6)$

Let ${\cal T}=\{\{j,3+j,6+j\},\{j,3+j,7+j\},\{j,2+j,6+j\}: j\in Z_{11}\}$, reducing the sums modulo $11$.  For each triple $T\in{\cal T}$, construct a $(3,1)$-GDD of type $8^3$ on the set $Z_8\times T$ with group set $\{Z_8\times \{l\}:l\in T\}$. Denote by ${\cal A}_T$ the set of its blocks. Let ${\cal B}_1=\bigcup_{T\in {\cal T}}{\cal A}_T$. ${\cal B}_2$ consists of $88\times35$ blocks which can be obtained by developing the following $35$ blocks by $(+1\ {\rm mod}\ 8,+1\ {\rm mod}\ 11)$:

\begin{tabular}{llllll}
$\{0_0,2_0,7_1\}$, & $\{0_0,2_0,2_2\}$, & $2\{0_0,2_0,3_2\}$,&
$\{0_0,3_0,6_2\}$, & $\{0_0,3_0,7_2\}$, & $\{0_0,4_0,5_2\}$,\\
$\{0_0,2_0,2_5\}$, & $\{0_0,2_0,0_5\}$, & $\{0_0,3_0,7_5\}$,&
$\{0_0,3_0,6_5\}$, & $\{0_0,3_0,5_5\}$, & $\{0_0,3_0,4_5\}$,\\
$\{0_0,4_0,5_5\}$, & $\{0_0,4_0,7_5\}$, & $\{0_0,1_1,3_2\}$,&
$\{0_0,2_1,5_2\}$, & $\{0_0,3_1,7_2\}$, & $\{0_0,6_1,5_2\}$,\\
$\{0_0,7_1,0_2\}$, & $\{0_0,1_1,2_2\}$, & $\{0_0,2_1,4_2\}$,&
$\{0_0,3_1,6_2\}$, & $\{0_0,4_1,0_2\}$, & $\{0_0,5_1,2_2\}$,\\
$\{0_0,6_1,4_2\}$, & $\{0_0,7_1,6_2\}$, & $\{0_0,1_1,7_2\}$,&
$\{0_0,2_1,7_2\}$, & $\{0_0,2_1,6_2\}$, & $\{0_0,1_1,5_2\}$,\\
$\{0_0,5_1,4_2\}$, & $\{0_0,3_1,1_2\}$, & $\{0_0,3_1,0_2\}$,&
$\{0_0,4_1,2_2\}$.
\end{tabular}

\noindent Then $(2\cdot{\cal B}_1)\cup{\cal B}_2$ forms a $2$-BSEC$(8,11,3,6)$.

\section{Appendix}

We here construct a $2$-BSEC$(r,c,3,2)$ for each $r\in\{11,14,17,20\}$ and $c\in\{6,10\}$.

When $c=6$, let $\cal A$ consist of $2r$ blocks obtained by developing the following $2$ blocks by $(+1\ {\rm mod}\ r,-)$: $\{0_0,0_2,0_4\}$ and $\{0_1,0_3,0_5\}$.
When $c=10$, let $\cal A$ consist of $10r$ blocks obtained by developing the following block by $(+1\ {\rm mod}\ r,$ $+1\ {\rm mod}\ 10)$: $\{0_0,0_2,0_6\}$.

$\bullet$ $r=11$

Let ${\cal T}_1=\{\{j,1+j,4+j\},\{j,2+j,5+j\}:j\in Z_{11}\}$, reducing the sums modulo $11$. For each triple $T\in{\cal T}_1$, construct a $(3,1)$-HGDD of type $(3,2^{c/2})$ on the set $T\times Z_c$ with group set $\{\{l\}\times Z_c:l\in T\}$ and hole set $\{T\times\{l,c/2+l\}:0\leq l< c/2\}$. Denote by ${\cal A}_T$ the set of its blocks. Let ${\cal B}_1=\bigcup_{T\in {\cal T}_1}{\cal A}_T$.

Let ${\cal T}_2=\{\{j,2+j,6+j\}:j\in Z_{11}\}$, reducing the sums modulo $11$. For each triple $T\in{\cal T}_2$, construct a $(3,1)$-MGDD of type $c^3$ on the set $T\times Z_c$ with group set $\{\{l\}\times Z_c:l\in T\}$ and hole set $\{T\times\{l\}:l\in Z_c\}$. Denote by ${\cal A}_T$ the set of its blocks. Let ${\cal B}_2=\bigcup_{T\in {\cal T}_2}{\cal A}_T$.

Let ${\cal B}_3$ consist of $44c$ blocks obtained by developing the following blocks by $(+1\ {\rm mod}\ 11,$ $+1\ {\rm mod}\ c)$:

\begin{tabular}{llll}
$\{0_0,5_{c/2},10_{c/2}\}$, & $\{0_0,5_0,8_{c/2}\}$,&
$\{0_0,4_0,8_0\}$, & $\{0_0,1_{c/2},4_{c/2}\}$.
\end{tabular}

When $c=6$, let ${\cal B}_4$ consist of $198$ blocks obtained by developing the following blocks by $(+1\ {\rm mod}\ 11,$ $+1\ {\rm mod}\ 6)$:

\begin{tabular}{lll}
$\{0_0,1_1,2_0\}$, & $\{0_0,0_2,1_4\}$, &
$\{0_0,0_3,2_3\}$.
\end{tabular}

\noindent When $c=10$, let ${\cal B}_4$ consist of $550$ blocks obtained by developing the following blocks by $(+1\ {\rm mod}\ 11,$ $+1\ {\rm mod}\ 10)$:

\begin{tabular}{lllll}
$\{0_0,0_2,0_7\}$, & $\{0_0,0_3,1_4\}$,&
$\{0_0,1_9,2_5\}$, & $\{0_0,1_3,2_0\}$, &
$\{0_0,1_2,2_0\}$.
\end{tabular}

Then ${\cal A}\cup (\bigcup_{i=1}^4{\cal B}_i)$ forms a $2$-BSEC$(11,c,3,2)$.

$\bullet$ $r=14$

Let ${\cal T}_1=\{\{j,1+j,6+j\},\{j,2+j,4+j\}:j\in Z_{14}\}$, reducing the sums modulo $14$. For each triple $T\in{\cal T}_1$, construct a $(3,1)$-HGDD of type $(3,2^{c/2})$ on the set $T\times Z_c$ with group set $\{\{l\}\times Z_c:l\in T\}$ and hole set $\{T\times\{l,c/2+l\}:0\leq l< c/2\}$. Denote by ${\cal A}_T$ the set of its blocks. Let ${\cal B}_1=\bigcup_{T\in {\cal T}_1}{\cal A}_T$.

Let ${\cal T}_2=\{\{j,3+j,7+j\},\{j,3+j,8+j\}:j\in Z_{14}\}$, reducing the sums modulo $14$. For each triple $T\in{\cal T}_2$, construct a $(3,1)$-GDD of type $c^3$ on the set $T\times Z_c$ with group set $\{\{l\}\times Z_c:l\in T\}$. Denote by ${\cal A}_T$ the set of its blocks. Let ${\cal B}_2=\bigcup_{T\in {\cal T}_2}{\cal A}_T$.

Let ${\cal B}_3$ consist of $42c$ blocks obtained by developing the following blocks by $(+1\ {\rm mod}\ 14,$ $+1\ {\rm mod}\ c)$:

\begin{tabular}{lll}
$\{0_0,1_{c/2},5_0\}$, & $\{0_0,1_{c/2},6_0\}$,&
$\{0_0,2_{c/2},6_{c/2}\}$.
\end{tabular}

When $c=6$, let ${\cal B}_4$ consist of $252$ blocks obtained by developing the following blocks by $(+1\ {\rm mod}\ 14,$ $+1\ {\rm mod}\ 6)$:

\begin{tabular}{lll}
$\{0_0,1_1,2_0\}$, & $\{0_0,0_2,1_4\}$,& $\{0_0,0_3,2_3\}$.
\end{tabular}

\noindent When $c=10$, let ${\cal B}_4$ consist of $700$ blocks obtained by developing the following blocks by $(+1\ {\rm mod}\ 14,$ $+1\ {\rm mod}\ 10)$:

\begin{tabular}{lllll}
$\{0_0,0_3,1_{c/2-1}\}$, & $\{0_0,1_{c-1},2_{c/2}\}$,&
$\{0_0,1_{c/2-2},2_0\}$, & $\{0_0,1_{c/2-3},2_0\}$,&
$\{0_0,0_2,0_7\}$.
\end{tabular}

Then ${\cal A}\cup (\bigcup_{i=1}^4{\cal B}_i)$ forms a $2$-BSEC$(14,c,3,2)$.

$\bullet$ $r=17$

Let ${\cal T}_1=\{\{j,1+j,3+j\},\{j,2+j,7+j\}:j\in Z_{17}\}$, reducing the sums modulo $17$. For each triple $T\in{\cal T}_1$, construct a $(3,1)$-HGDD of type $(3,2^{c/2})$ on the set $T\times Z_c$ with group set $\{\{l\}\times Z_c:l\in T\}$ and hole set $\{T\times\{l,c/2+l\}:0\leq l< c/2\}$. Denote by ${\cal A}_T$ the set of its blocks. Let ${\cal B}_1=\bigcup_{T\in {\cal T}_1}{\cal A}_T$.

Let ${\cal T}_2=\{\{j,4+j,9+j\},\{j,3+j,9+j\}:j\in Z_{17}\}$, reducing the sums modulo $17$. For each triple $T\in{\cal T}_2$, construct a $(3,1)$-MGDD of type $c^3$ on the set $T\times Z_c$ with group set $\{\{l\}\times Z_c:l\in T\}$ and hole set $\{T\times\{l\}:l\in Z_c\}$. Denote by ${\cal A}_T$ the set of its blocks. Let ${\cal B}_2=\bigcup_{T\in {\cal T}_2}{\cal A}_T$.

Let ${\cal T}_3=\{\{j,4+j,10+j\}:j\in Z_{17}\}$, reducing the sums modulo $17$. For each triple $T\in{\cal T}_3$, construct a $(3,1)$-GDD of type $c^3$ on the set $T\times Z_c$ with group set $\{\{l\}\times Z_c:l\in T\}$. Denote by ${\cal A}_T$ the set of its blocks. Let ${\cal B}_3=\bigcup_{T\in {\cal T}_3}{\cal A}_T$.

Let ${\cal B}_4$ consist of $85c$ blocks obtained by developing the following blocks by $(+1\ {\rm mod}\ 17,$ $+1\ {\rm mod}\ c)$:

\begin{tabular}{lllll}
$\{0_0,1_{c/2},8_0\}$, & $\{0_0,1_{c/2},6_0\}$,&
$\{0_0,3_0,7_0\}$, & $\{0_0,2_{c/2},5_0\}$,&
$\{0_0,3_0,8_0\}$.
\end{tabular}

When $c=6$, let ${\cal B}_5$ consist of $306$ blocks obtained by developing the following blocks by $(+1\ {\rm mod}\ 17,$ $+1\ {\rm mod}\ 6)$:

\begin{tabular}{lll}
$\{0_0,1_1,2_0\}$, & $\{0_0,0_2,1_4\}$,& $\{0_0,0_3,2_3\}$.
\end{tabular}

\noindent When $c=10$, let ${\cal B}_5$ consist of $850$ blocks obtained by developing the following blocks by $(+1\ {\rm mod}\ 17,$ $+1\ {\rm mod}\ 10)$:

\begin{tabular}{lllll}
$\{0_0,0_3,1_4\}$, & $\{0_0,1_9,2_5\}$,& $\{0_0,1_3,2_0\}$, & $\{0_0,1_2,2_0\}$,&
$\{0_0,0_2,0_7\}$.
\end{tabular}

Then ${\cal A}\cup (\bigcup_{i=1}^5{\cal B}_i)$ forms a $2$-BSEC$(17,c,3,2)$.

$\bullet$ $r=20$

Let ${\cal T}_1=\{\{j,1+j,2+j\},\{j,5+j,12+j\}:j\in Z_{20}\}$, reducing the sums modulo $20$. For each triple $T\in{\cal T}_1$, construct a $(3,1)$-HGDD of type $(3,2^{c/2})$ on the set $T\times Z_c$ with group set $\{\{l\}\times Z_c:l\in T\}$ and hole set $\{T\times\{l,c/2+l\}:0\leq l< c/2\}$. Denote by ${\cal A}_T$ the set of its blocks. Let ${\cal B}_1=\bigcup_{T\in {\cal T}_1}{\cal A}_T$.

Let ${\cal T}_2=\{\{j,2+j,10+j\},\{j,3+j,6+j\}:j\in Z_{20}\}$, reducing the sums modulo $20$. For each triple $T\in{\cal T}_2$, construct a $(3,1)$-MGDD of type $c^3$ on the set $T\times Z_c$ with group set $\{\{l\}\times Z_c:l\in T\}$ and hole set $\{T\times\{l\}:l\in Z_c\}$. Denote by ${\cal A}_T$ the set of its blocks. Let ${\cal B}_2=\bigcup_{T\in {\cal T}_2}{\cal A}_T$.

Let ${\cal T}_3=\{\{j,4+j,11+j\},\{j,5+j,11+j\}:j\in Z_{20}\}$, reducing the sums modulo $20$. For each triple $T\in{\cal T}_3$, construct a $(3,1)$-GDD of type $c^3$ on the set $T\times Z_c$ with group set $\{\{l\}\times Z_c:l\in T\}$. Denote by ${\cal A}_T$ the set of its blocks. Let ${\cal B}_3=\bigcup_{T\in {\cal T}_3}{\cal A}_T$.

Let ${\cal B}_4$ consist of $100c$ blocks obtained by developing the following blocks by $(+1\ {\rm mod}\ 20,$ $+1\ {\rm mod}\ c)$:

\begin{tabular}{lllll}
$\{0_0,1_{c/2},2_0\}$, & $\{0_0,2_0,5_0\}$,&
$\{0_0,4_0,10_0\}$, & $\{0_0,2_{c/2},7_0\}$,&
$\{0_0,3_0,7_{c/2}\}$.
\end{tabular}

When $c=6$, let ${\cal B}_5$ consist of $360$ blocks obtained by developing the following blocks by $(+1\ {\rm mod}\ 20,$ $+1\ {\rm mod}\ 6)$:

\begin{tabular}{lll}
$\{0_0,0_4,4_5\}$, & $\{0_0,0_3,8_3\}$,& $\{0_0,4_2,8_0\}$.
\end{tabular}

\noindent When $c=10$, let ${\cal B}_5$ consist of $1000$ blocks obtained by developing the following blocks by $(+1\ {\rm mod}\ 20,$ $+1\ {\rm mod}\ 10)$:

\begin{tabular}{lllll}
$\{0_0,0_3,4_{c/2-1}\}$, & $\{0_0,4_{c-1},8_{c/2}\}$,&
$\{0_0,4_{c/2-2},8_0\}$, & $\{0_0,4_{c/2-3},8_0\}$,&
$\{0_0,0_2,0_7\}$.
\end{tabular}

Then ${\cal A}\cup (\bigcup_{i=1}^5{\cal B}_i)$ forms a $2$-BSEC$(20,c,3,2)$.

\section{Appendix}

We here construct a $2$-BSEC$(r,c,3,2)$ for $(r,c)\in\{(3,10),(3,14),(3,16),(3,20),(3,22),(5,6),\linebreak(5,10)\}$. We found them by computer search. Conveniently, all elements of the required $2$-BSEC$(r,c,3,2)$ are defined on $Z_{rc}$, and arranged as follows
\begin{center}
$
         \begin{array}{ccccc}
           0 & r & 2r & \cdots & (c-1)r \\
           1 & r+1 & 2r+1 & \cdots & (c-1)r +1 \\
           \vdots&\vdots&\vdots&\vdots&\vdots\\
           r-1 & 2r-1 & 3r-1 & \cdots & rc-1 \\
         \end{array}
$\end{center}
It is required that every pair of consecutive points in each row and each column do not appear in any block. Let $\alpha$ be a permutation on $Z_{rc}$ and $G$ be the group generated by $\alpha$ (or let $\alpha$ and $\beta$ be two permutations on $Z_{rc}$ and $G$ be the group generated by $\alpha$ and $\beta$). Only initial blocks are listed below. All
other blocks are obtained by developing these initial blocks under the
action of $G$.

$\bullet$ $(r,c)=(3,10)$

$\alpha=(0\ 6\ 12\ 18\ 24)(1\ 7\ 13\ 19\ 25)\cdots(5\ 11\ 17\ 23\ 29)$

\begin{tabular}{llllll}
$\{0,20,25\}$&
$\{0,21,26\}$&
$\{0,22,26\}$&
$\{1,3,8\}$&
$\{1,5,9\}$&
$\{1,7,13\}$\\
$\{1,8,14\}$&
$\{1,9,15\}$&
$\{1,10,15\}$&
$\{1,10,19\}$&
$\{1,11,20\}$&
$\{1,11,20\}$\\
$\{1,14,21\}$&
$\{1,16,22\}$&
$\{1,16,23\}$&
$\{1,17,21\}$&
$\{1,17,23\}$&
$\{2,8,17\}$\\
$\{2,9,16\}$&
$\{2,10,17\}$&
$\{2,10,20\}$&
$\{2,11,16\}$&
$\{2,15,21\}$&
$\{2,15,22\}$\\
$\{3,11,21\}$&
$\{3,11,7\}$&
$\{0,4,6\}$&
$\{0,4,6\}$&
$\{0,5,7\}$&
$\{0,5,15\}$\\
$\{0,7,18\}$&
$\{0,8,9\}$&
$\{0,8,20\}$&
$\{0,9,23\}$&
$\{0,10,14\}$&
$\{0,10,21\}$\\
$\{0,11,22\}$&
$\{0,11,23\}$&
$\{0,12,25\}$&
$\{0,13,15\}$&
$\{0,14,16\}$&
$\{0,16,29\}$\\
$\{0,17,19\}$&
$\{0,17,29\}$&
$\{1,26,27\}$&
$\{2,4,21\}$&
$\{3,15,28\}$&
$\{3,17,22\}$\\
$\{4,10,22\}$&
$\{4,17,23\}$
\end{tabular}

$\bullet$ $(r,c)=(3,14)$

$\alpha=(0\ 6\ 12\ 18\ 24\ 30\ 36)(1\ 7\ 13\ 19\ 25\ 31\ 37)\cdots(5\ 11\ 17\ 23\ 29\ 35\ 41)$

\begin{tabular}{llllll}
$\{2,14,26\}$&
$\{2,15,23\}$&
$\{2,15,26\}$&
$\{2,16,27\}$&
$\{2,16,28\}$&
$\{2,17,27\}$\\
$\{2,17,29\}$&
$\{2,21,28\}$&
$\{2,21,29\}$&
$\{2,22,33\}$&
$\{2,23,35\}$&
$\{3,15,27\}$\\
$\{3,16,18\}$&
$\{0,4,6\}$&
$\{0,4,8\}$&
$\{0,5,6\}$&
$\{0,5,7\}$&
$\{0,7,9\}$\\
$\{0,8,9\}$&
$\{0,10,14\}$&
$\{0,10,15\}$&
$\{0,11,12\}$&
$\{0,11,13\}$&
$\{0,12,25\}$\\
$\{0,14,15\}$&
$\{0,16,24\}$&
$\{0,16,28\}$&
$\{0,17,34\}$&
$\{0,17,35\}$&
$\{0,18,37\}$\\
$\{0,19,21\}$&
$\{0,20,22\}$&
$\{0,20,22\}$&
$\{0,21,35\}$&
$\{0,23,28\}$&
$\{0,23,29\}$\\
$\{0,25,29\}$&
$\{0,26,31\}$&
$\{0,26,32\}$&
$\{0,27,31\}$&
$\{0,32,37\}$&
$\{3,16,29\}$\\
$\{3,17,21\}$&
$\{3,22,28\}$&
$\{3,22,40\}$&
$\{3,23,7\}$&
$\{0,33,38\}$&
$\{0,33,38\}$\\
$\{1,5,9\}$&
$\{1,7,13\}$&
$\{1,8,14\}$&
$\{1,8,15\}$&
$\{1,9,15\}$&
$\{1,10,16\}$\\
$\{1,10,17\}$&
$\{1,11,25\}$&
$\{1,11,27\}$&
$\{1,13,34\}$&
$\{1,14,25\}$&
$\{1,16,35\}$\\
$\{1,20,28\}$&
$\{1,20,29\}$&
$\{1,21,27\}$&
$\{1,21,28\}$&
$\{1,22,32\}$&
$\{1,23,33\}$\\
$\{1,23,34\}$&
$\{1,26,33\}$&
$\{1,26,35\}$&
$\{2,10,34\}$&
$\{2,22,35\}$&
$\{3,23,28\}$\\
$\{4,11,35\}$&
$\{4,23,29\}$
\end{tabular}

$\bullet$ $(r,c)=(3,16)$

$\alpha=(0\ 6\ 12\ \cdots\ 42)(1\ 7\ 13\ \cdots\ 43)\cdots(5\ 11\ 17\ \cdots\ 47)$

\begin{tabular}{llllll}
$\{0,4,24\}$&
$\{0,4,28\}$&
$\{1,5,25\}$&
$\{1,5,29\}$&
$\{2,3,26\}$&
$\{2,3,27\}$\\
$\{0,35,42\}$&
$\{0,5,36\}$&
$\{0,5,47\}$&
$\{0,6,18\}$&
$\{0,7,29\}$&
$\{0,7,33\}$\\
$\{0,8,43\}$&
$\{0,8,18\}$&
$\{0,9,29\}$&
$\{0,9,26\}$&
$\{0,10,20\}$&
$\{0,10,39\}$\\
$\{0,11,13\}$&
$\{0,11,40\}$&
$\{0,13,15\}$&
$\{0,14,47\}$&
$\{0,14,46\}$&
$\{0,15,46\}$\\
$\{0,16,39\}$&
$\{0,16,25\}$&
$\{0,17,25\}$&
$\{0,19,31\}$&
$\{0,19,41\}$&
$\{0,20,22\}$\\
$\{0,21,32\}$&
$\{0,21,26\}$&
$\{0,22,37\}$&
$\{0,23,35\}$&
$\{0,23,37\}$&
$\{0,27,33\}$\\
$\{0,27,31\}$&
$\{0,32,34\}$&
$\{0,34,43\}$&
$\{0,38,44\}$&
$\{0,40,44\}$&
$\{1,3,7\}$\\
$\{1,7,15\}$&
$\{1,8,20\}$&
$\{1,8,21\}$&
$\{1,9,15\}$&
$\{1,10,31\}$&
$\{1,10,37\}$\\
$\{1,11,22\}$&
$\{1,11,38\}$&
$\{1,14,26\}$&
$\{1,16,31\}$&
$\{1,16,32\}$&
$\{1,17,38\}$\\
$\{1,17,47\}$&
$\{1,20,39\}$&
$\{1,21,26\}$&
$\{1,27,44\}$&
$\{1,28,32\}$&
$\{1,33,41\}$\\
$\{1,33,44\}$&
$\{1,35,39\}$&
$\{2,8,15\}$&
$\{2,9,16\}$&
$\{2,10,16\}$&
$\{2,10,17\}$\\
$\{2,11,17\}$&
$\{2,11,20\}$&
$\{2,20,40\}$&
$\{2,21,28\}$&
$\{2,22,29\}$&
$\{2,23,35\}$\\
$\{2,28,41\}$&
$\{3,11,15\}$&
$\{3,15,28\}$&
$\{3,16,21\}$&
$\{3,17,33\}$&
$\{3,17,34\}$\\
$\{3,22,41\}$&
$\{3,23,41\}$&
$\{3,29,40\}$&
$\{3,29,46\}$&
$\{3,35,40\}$&
$\{4,10,22\}$\\
$\{4,16,29\}$&
$\{4,22,47\}$
\end{tabular}

$\bullet$ $(r,c)=(3,20)$

$\alpha=(0\ 6\ 12\ \cdots\ 54)(1\ 7\ 13\ \cdots\ 55)\cdots(5\ 11\ 17\ \cdots\ 59)$

\begin{tabular}{llllll}
$\{0,4,30\}$&
$\{0,4,34\}$&
$\{1,5,31\}$&
$\{1,5,35\}$&
$\{2,3,32\}$&
$\{2,3,33\}$\\
$\{0,35,42\}$&
$\{0,5,48\}$&
$\{0,5,41\}$&
$\{0,6,31\}$&
$\{0,6,19\}$&
$\{0,7,19\}$\\
$\{0,7,23\}$&
$\{0,8,45\}$&
$\{0,8,28\}$&
$\{0,9,56\}$&
$\{0,9,28\}$&
$\{0,10,37\}$\\
$\{0,10,52\}$&
$\{0,11,25\}$&
$\{0,11,12\}$&
$\{0,13,35\}$&
$\{0,14,53\}$&
$\{0,14,23\}$\\
$\{0,15,22\}$&
$\{0,15,58\}$&
$\{0,16,39\}$&
$\{0,16,36\}$&
$\{0,17,49\}$&
$\{0,18,55\}$\\
$\{0,20,49\}$&
$\{0,20,56\}$&
$\{0,21,55\}$&
$\{0,21,52\}$&
$\{0,22,24\}$&
$\{0,26,38\}$\\
$\{0,26,38\}$&
$\{0,27,45\}$&
$\{0,27,32\}$&
$\{0,29,50\}$&
$\{0,29,46\}$&
$\{0,31,50\}$\\
$\{0,32,39\}$&
$\{0,33,43\}$&
$\{0,33,40\}$&
$\{0,41,43\}$&
$\{0,44,51\}$&
$\{0,44,46\}$\\
$\{0,47,59\}$&
$\{0,47,51\}$&
$\{1,3,11\}$&
$\{1,3,26\}$&
$\{1,7,38\}$&
$\{1,7,15\}$\\
$\{1,8,25\}$&
$\{1,8,22\}$&
$\{1,9,14\}$&
$\{1,10,37\}$&
$\{1,10,26\}$&
$\{1,11,19\}$\\
$\{1,13,51\}$&
$\{1,14,33\}$&
$\{1,15,57\}$&
$\{1,16,22\}$&
$\{1,16,39\}$&
$\{1,17,19\}$\\
$\{1,20,41\}$&
$\{1,21,38\}$&
$\{1,21,29\}$&
$\{1,23,56\}$&
$\{1,27,40\}$&
$\{1,28,46\}$\\
$\{1,28,52\}$&
$\{1,33,52\}$&
$\{1,40,53\}$&
$\{1,41,47\}$&
$\{1,44,50\}$&
$\{1,45,56\}$\\
$\{1,45,57\}$&
$\{1,46,50\}$&
$\{2,4,38\}$&
$\{2,8,40\}$&
$\{2,10,15\}$&
$\{2,10,16\}$\\
$\{2,11,17\}$&
$\{2,17,21\}$&
$\{2,20,40\}$&
$\{2,20,45\}$&
$\{2,23,28\}$&
$\{2,27,47\}$\\
$\{2,29,51\}$&
$\{2,34,53\}$&
$\{2,35,52\}$&
$\{2,35,53\}$&
$\{2,46,58\}$&
$\{2,47,52\}$\\
$\{3,9,15\}$&
$\{3,16,23\}$&
$\{3,17,27\}$&
$\{3,17,27\}$&
$\{3,28,47\}$&
$\{3,28,52\}$\\
$\{3,29,52\}$&
$\{3,29,58\}$&
$\{3,34,41\}$&
$\{3,35,46\}$&
$\{3,35,47\}$&
$\{4,16,41\}$\\
$\{4,17,35\}$&
$\{4,29,53\}$
\end{tabular}

$\bullet$ $(r,c)=(3,22)$

$\alpha=(0\ 6\ 12\ \cdots\ 60)(1\ 7\ 13\ \cdots\ 61)\cdots(5\ 11\ 17\ \cdots\ 65)$

\begin{tabular}{llllll}
$\{0,35,41\}$&
$\{0,37,43\}$&
$\{0,39,44\}$&
$\{0,39,44\}$&
$\{0,40,46\}$&
$\{0,45,52\}$\\
$\{0,45,56\}$&
$\{0,46,53\}$&
$\{0,47,58\}$&
$\{0,47,58\}$&
$\{0,50,55\}$&
$\{0,50,57\}$\\
$\{2,15,28\}$&
$\{2,16,32\}$&
$\{2,16,33\}$&
$\{2,17,33\}$&
$\{2,17,34\}$&
$\{2,20,39\}$\\
$\{2,21,35\}$&
$\{2,22,40\}$&
$\{2,22,40\}$&
$\{2,23,39\}$&
$\{2,26,13\}$&
$\{1,9,15\}$\\
$\{1,10,19\}$&
$\{1,10,20\}$&
$\{1,11,20\}$&
$\{1,11,21\}$&
$\{1,13,29\}$&
$\{1,15,22\}$\\
$\{1,16,26\}$&
$\{1,16,32\}$&
$\{1,17,32\}$&
$\{1,19,43\}$&
$\{1,21,38\}$&
$\{1,22,33\}$\\
$\{1,23,31\}$&
$\{1,23,35\}$&
$\{1,26,37\}$&
$\{1,27,45\}$&
$\{1,27,50\}$&
$\{1,29,46\}$\\
$\{1,33,47\}$&
$\{1,34,46\}$&
$\{1,34,57\}$&
$\{0,4,6\}$&
$\{0,4,8\}$&
$\{0,5,7\}$\\
\end{tabular}

\begin{tabular}{llllll}
$\{0,5,60\}$&
$\{0,7,9\}$&
$\{0,8,9\}$&
$\{0,10,12\}$&
$\{0,10,54\}$&
$\{0,11,13\}$\\
$\{0,13,17\}$&
$\{0,14,15\}$&
$\{0,14,65\}$&
$\{0,15,19\}$&
$\{2,35,47\}$&
$\{2,41,48\}$\\
$\{0,18,36\}$&
$\{0,19,23\}$&
$\{0,20,22\}$&
$\{0,26,51\}$&
$\{0,51,59\}$&
$\{0,23,52\}$\\
$\{0,24,49\}$&
$\{0,24,49\}$&
$\{0,26,30\}$&
$\{0,28,57\}$&
$\{0,53,61\}$&
$\{0,55,62\}$\\
$\{0,56,61\}$&
$\{1,9,41\}$&
$\{1,13,63\}$&
$\{1,14,44\}$&
$\{1,44,52\}$&
$\{1,45,58\}$\\
$\{2,8,26\}$&
$\{2,10,14\}$&
$\{2,11,56\}$&
$\{2,15,41\}$&
$\{2,46,59\}$&
$\{3,9,33\}$\\
$\{0,28,43\}$&
$\{0,29,33\}$&
$\{0,29,33\}$&
$\{0,31,37\}$&
$\{0,31,38\}$&
$\{3,11,34\}$\\
$\{3,15,35\}$&
$\{3,15,41\}$&
$\{3,22,47\}$&
$\{0,16,21\}$&
$\{0,16,27\}$&
$\{0,17,35\}$\\
$\{0,21,41\}$&
$\{0,27,65\}$&
$\{0,32,34\}$&
$\{0,32,38\}$&
$\{0,34,40\}$&
$\{1,3,40\}$\\
$\{1,28,51\}$&
$\{1,28,53\}$&
$\{1,35,41\}$&
$\{1,38,47\}$&
$\{1,39,50\}$&
$\{1,39,57\}$\\
$\{1,40,53\}$&
$\{2,9,45\}$&
$\{2,27,46\}$&
$\{2,28,47\}$&
$\{2,29,34\}$&
$\{2,29,51\}$\\
$\{3,27,52\}$&
$\{3,28,59\}$&
$\{3,34,64\}$&
$\{4,11,41\}$&
$\{4,16,40\}$&
$\{4,23,46\}$\\
$\{4,35,65\}$&
$\{5,23,47\}$
\end{tabular}

$\bullet$ $(r,c)=(5,6)$

$\alpha=(0\ 1\ 2\ 3\ 4) (5\ 6\ 7\ 8\ 9)\cdots(25\ 26\ 27\ 28\ 29)$

\begin{tabular}{llllll}
$\{0,2,12\}$&
$\{0,2,6\}$&
$\{0,6,19\}$&
$\{0,7,19\}$&
$\{0,7,23\}$&
$\{0,8,15\}$\\
$\{0,8,28\}$&
$\{0,9,10\}$&
$\{0,11,26\}$&
$\{0,11,28\}$&
$\{0,12,22\}$&
$\{0,13,26\}$\\
$\{0,13,23\}$&
$\{0,14,22\}$&
$\{0,14,27\}$&
$\{0,15,18\}$&
$\{0,16,20\}$&
$\{0,16,24\}$\\
$\{0,17,20\}$&
$\{0,17,27\}$&
$\{0,18,29\}$&
$\{0,21,24\}$&
$\{0,21,29\}$&
$\{5,7,15\}$\\
$\{5,7,22\}$&
$\{5,11,24\}$&
$\{5,12,29\}$&
$\{5,12,14\}$&
$\{5,13,15\}$&
$\{5,13,16\}$\\
$\{5,14,16\}$&
$\{5,19,20\}$&
$\{5,19,21\}$&
$\{5,22,26\}$&
$\{5,23,27\}$&
$\{5,23,29\}$\\
$\{5,24,27\}$&
$\{5,25,28\}$&
$\{5,26,28\}$&
$\{10,12,21\}$&
$\{10,16,22\}$&
$\{10,16,25\}$\\
$\{10,18,26\}$&
$\{10,19,26\}$&
$\{10,19,29\}$&
$\{10,21,24\}$&
$\{10,22,29\}$&
$\{15,17,29\}$\\
$\{15,22,28\}$&
$\{15,24,26\}$
\end{tabular}

$\bullet$ $(r,c)=(5,10)$

$\alpha=(0\ 1\ 2\ 3\ 4) (5\ 6\ 7\ 8\ 9)\cdots(45\ 46\ 47\ 48\ 49)$

$\beta=(0\ 10\ 20\ 30\ 40)(5\ 15\ 25\ 35\ 45)\cdots(4\ 14\ 24\ 34\ 44)(9\ 19\ 29\ 39\ 49)$

\begin{tabular}{llllll}
$\{0,2,42\}$&
$\{0,2,18\}$&
$\{0,6,35\}$&
$\{0,6,20\}$&
$\{0,7,25\}$&
$\{0,7,29\}$\\
$\{0,8,15\}$&
$\{0,8,46\}$&
$\{0,9,32\}$&
$\{0,9,28\}$&
$\{0,10,42\}$&
$\{0,11,46\}$\\
$\{0,11,28\}$&
$\{0,12,37\}$&
$\{0,12,33\}$&
$\{0,14,47\}$&
$\{0,14,48\}$&
$\{0,15,22\}$\\
$\{0,16,39\}$&
$\{0,17,36\}$&
$\{0,18,24\}$&
$\{0,19,31\}$&
$\{0,19,20\}$&
$\{0,21,47\}$\\
$\{0,26,37\}$&
$\{0,27,29\}$&
$\{5,7,18\}$&
$\{5,15,45\}$&
$\{5,18,46\}$&
$\{5,19,39\}$
\end{tabular}


\makeatletter
\renewcommand\@biblabel[1]{}
\renewenvironment{thebibliography}[1]
     {\section*{\refname}%
      \@mkboth{\MakeUppercase\refname}{\MakeUppercase\refname}%
      \list{\@biblabel{\@arabic\c@enumiv}}%
           {\settowidth\labelwidth{\@biblabel{#1}}%
            \leftmargin\labelwidth
            \advance\leftmargin by 2em%
            \itemindent -2em%
            \@openbib@code
            \usecounter{enumiv}%
            \let\p@enumiv\@empty
            \renewcommand\theenumiv{\@arabic\c@enumiv}}%
      \sloppy
      \clubpenalty4000
      \@clubpenalty \clubpenalty
      \widowpenalty4000%
      \sfcode`\.\@m}
     {\def\@noitemerr
       {\@latex@warning{Empty `thebibliography' environment}}%
      \endlist}
\makeatother


\end{document}